\documentclass[8pt]{article}

\usepackage{amsmath,amssymb,amsfonts,amsthm,amscd,shadow,fancybox,epic}
\usepackage{youngtab}
\usepackage{pstricks,pst-node,pst-plot}

\theoremstyle{plain}

\newtheorem{lem}[subsection]{Lemma}
\newtheorem{prop}[subsection]{Proposition}

\newtheorem{ex}[subsection]{Example}

\theoremstyle{definition}
\newtheorem{rem}[subsection]{Remark}

\newcommand{\gl}{\lambda}
\newcommand{\gm}{\mu}
\newcommand{\gn}{\nu}
\newcommand{\gp}{\pi}
\newcommand{\gt}{\theta}
\newcommand{\ga}{\alpha}
\newcommand{\gb}{\beta}
\newcommand{\gc}{\gamma}

\newcommand{\ft}{\mathfrak{t}}
\newcommand{\bC}{\mathbb{C}}
\newcommand{\cE}{\mathcal{E}}
\newcommand{\cF}{\mathcal{F}}
\newcommand{\cD}{\mathcal{D}}
\newcommand{\cP}{\st}
\newcommand{\cM}{\mathcal{M}}
\newcommand{\cN}{\mathcal{N}}
\newcommand{\cX}{\mathcal{X}}
\newcommand{\cY}{\mathcal{Y}}
\newcommand{\cZ}{\mathcal{Z}}

\newcommand{\schub}{\mathfrak{S}}
\newcommand{\groth}{\mathcal{G}}
\newcommand{\bb}{\mathbb}
\newcommand{\supp}{\mathop{\rm Supp}\nolimits}
\newcommand{\pl}{\mathop{\rm pl}\nolimits}
\newcommand{\Span}{\mathop{\rm Span}\nolimits}
\newcommand{\Char}{\mathop{\rm Char}\nolimits}

\newcommand{\sstab}{\cP}

\newcommand{\ssvt}{\mathop{\rm SSVT}\nolimits}
\newcommand{\st}{\mathop{\rm SSYT}\nolimits}
\newcommand{\svt}{\mathop{\rm SVT}\nolimits}

\newcommand{\glb}{\mathop{\rm glb}\nolimits}
\newcommand{\tw}{\mathop{\rm twist}\nolimits}

\newcommand{\ol}{\overline}

\newcommand{\sK}{_{\text{\tiny{K}}}}
\newcommand{\sH}{_{\text{\tiny{H}}}}

\newcommand{\cO}{\mathcal{O}}
\newcommand{\ab}{_{\ga,\gb}}
\newcommand{\lm}{_{\gl,\gm}}
\newcommand{\tp}{^{\,\rm top}}

\title{Schubert Classes in the
Equivariant K-Theory and Equivariant Cohomology\\ of the
Grassmannian}

\author{Victor Kreiman}

\begin{document}
\maketitle

\begin{abstract}
We give positive formulas for the restriction of a Schubert Class
to a $T$-fixed point in the equivariant K-theory and equivariant
cohomology of the Grassmannian.  Our formulas rely on a result of
Kodiyalam-Raghavan and Kreiman-Lakshmibai, which gives an
equivariant Gr\"obner degeneration of a Schubert variety in the
neighborhood of a $T$-fixed point of the Grassmannian.
\end{abstract}

\setcounter{tocdepth}{1} \tableofcontents

\section{Introduction}\label{s.intro}

The group $T$ of diagonal matrices in $GL_n(\bC)$ acts on the
Grassmannian $Gr_{d,n}$, with fixed point set indexed by
$I_{d,n}$, the $d$ element subsets of $\{1,\ldots,n\}$. For
$\gb\in I_{d,n}$, denote the corresponding $T$-fixed point by
$e_\gb$. The $T$-equivariant embedding $e_\gb\stackrel{i}{\to}
Gr_{d,n}$ induces restriction homomorphisms $i^*_K$ in
$T$-equivariant K-theory and $i^*_H$ in $T$-equivariant
cohomology:
\begin{align*}
K_T^*(Gr_{d,n})&\stackrel{i^*_K}{\rightarrow} K_T^*(e_\gb)\cong
R(T)=\bC[t_1^{\pm
1},\ldots,t_n^{\pm 1}]\\
H_T^*(Gr_{d,n})&\stackrel{i^*_H}{\rightarrow} H_T^*(e_\gb)\cong
\bC[\ft^*] =\bC[t_1,\ldots,t_n]
\end{align*}
where $R(T)$ is the representation ring of $T$ and $\ft$ is the
Lie algebra of $T$.  The image of an element $z$ of
$K_T^*(Gr_{d,n})$ or $H_T^*(Gr_{d,n})$ under restriction to
$e_\gb$ is denoted by $z|_{e_\gb}$.  The product maps
\begin{align*}
&K_T^*(Gr_{d,n})\to\prod\limits_{\gb\in I_{d,n}}K_T^*(e_\gb),
\ \ z\mapsto \prod\limits_{\gb\in I_{d,n}} z|_{e_\gb}\\
&H_T^*(Gr_{d,n})\to\prod\limits_{\gb\in I_{d,n}}H_T^*(e_\gb),\ \
z\mapsto \prod\limits_{\gb\in I_{d,n}} z|_{e_\gb}
\end{align*}
are both injective. Thus, an element of $K_T^*(Gr_{d,n})$ or
$H_T^*(Gr_{d,n})$ is determined by its restrictions to all
$e_\gb$, $\gb\in I_{d,n}$.

The Schubert varieties of the Grassmannian are in bijection with
the $T$-fixed points, and thus are also indexed by $I_{d,n}$. For
$\ga\in I_{d,n}$, denote the Schubert variety by $X_\ga$, and the
corresponding Schubert classes in $K_T^*(Gr_{d,n})$ and
$H_T^*(Gr_{d,n})$ by $[X_\ga]\sK$ and $[X_\ga]\sH$ respectively.
In this paper, we obtain formulas for $[X_\ga]\sK|_{e_\gb}$ and
$[X_\ga]\sH|_{e_\gb}$.

Various formulas already exist for $[X_\ga]\sK|_{e_\gb}$ and
$[X_\ga]\sH|_{e_\gb}$. Letting $\groth_\ga({\bf r},{\bf t})$ and
$\schub_\ga({\bf r},{\bf t})$ denote the double Grothendieck and
double Schubert polynomials \cite{La-Sc,La-Sc2} respectively for
$\ga$,
\begin{equation}\label{e.groth_schub}
[X_\ga]\sK|_{e_\gb}=\groth_\ga(\gb({\bf t}),{\bf t})\ \ \hbox{ and
}\ \ \ [X_\ga]\sH|_{e_\gb}=\schub_\ga(\gb({\bf t}),{\bf t}),
\end{equation}
where on the right sides of both equations we view $\ga$ and $\gb$
as Grassmannian permutations rather than $d$-tuples. Letting
$d_{\ga,\gb}^{\,\gc}$ and $c_{\ga,\gb}^{\,\gc}$ denote the linear
structure constants for the Schubert classes in the equivariant
K-theory and equivariant cohomology respectively of the
Grassmannian,
\begin{equation}
[X_\ga]\sK|_{e_\gb}=d_{\ga,\gb}^{\,\gb}\ \ \hbox{ and }\ \ \
[X_\ga]\sH|_{e_\gb}=c_{\ga,\gb}^{\,\gb}.
\end{equation}
Hence the various formulas for $\groth_\ga({\bf r},{\bf t})$ and
$\schub_\ga({\bf r},{\bf t})$ (see
\cite{Be,Be-Bi,Bi,Bu-Kr-Ta-Yo,Fo-Ki,Fo-Ki2,Kn-Mi-Yo2,Kn-Mi-Yo,Kn-Mi,La,La-Le-Th,La-Sc,La-Sc2,Le,Le-Ro-So,Ma},
for example) and  for $c_{\ga,\gb}^{\,\gc}$ (see
\cite{Kn-Ta,Mo-Sa}) can be used to compute $[X_\ga]\sK|_{e_\gb}$
and $[X_\ga]\sK|_{e_\gb}$. The main features of our formulas are
that they satisfy positivity conditions, they are obtained via a
Gr\"obner degeneration, and they are expressed in terms of
semistandard set-valued tableaux.

\subsection*{Posititivity}

Griffeth and Ram \cite{Gr-Ra} conjecture that the structure
constants $d_{\ga,\gb}^{\,\gc}$ for $G/B$, where $G$ is any
symmetrizable Kac-Moody group, can be expressed as\\
$(-1)^{l(\ga)+l(\gb)-l(\gc)}$ times a sum of products of terms of
the form $e^{\gt}-1$ or $e^{\gt}$, where $\gt$ is a positive root.
We prove and realize this positivity conjecture for
$d_{\ga,\gb}^{\,\gb}$, Grassmannian $G/P$. Our formula involves
only terms of the form $e^{\gt}-1$ (in our case, $t_b/t_a-1$,
$b>a$).

Graham \cite{Gr} proves that the structure constants
$c_{\ga,\gb}^{\,\gc}$ for $G/B$, where $G$ is any symmetrizable
Kac-Moody group, can be expressed as sums of products of positive
roots. We realize this condition for $c_{\ga,\gb}^{\,\gb}$,
Grassmannian $G/P$. Knutson and Tao \cite{Kn-Ta} realize this
positivity condition for all $c_{\ga,\gb}^{\,\gc}$, Grassmannian
$G/P$. Their formula, when restricted to structure constants of
the form $c_{\ga,\gb}^{\,\gb}$, is expressed in terms of different
combinatorial objects than ours, and also expresses the quantity
$c_{\ga,\gb}^{\,\gb}$ in terms of different sums of monomials in
the positive roots.

\subsection*{Gr\"obner Degeneration}

Our proof relies on a result of Kodiyalam and Raghavan
\cite{Ko-Ra}, Kreiman and Lakshmibai \cite{Kr-La}, and Kreiman
\cite{Kr}, which gives an equivariant Gr\"obner degeneration of a
local neighborhood of $X_\ga$ centered at $e_\gb$ to a reduced
union $W\ab$ of coordinate subspaces $W_1,\ldots,W_q$ of an affine
space whose coordinates are characters of $T$. Our strategy is to
use this result and the inclusion-exclusion principle to deduce
that
\begin{align*}
[X_\ga]\sK|_{e_\gb}=[W\ab]\sK&=\sum\limits_{j=1}^q(-1)^{j+1}\sum\limits_{1\leq
i_1<\cdots<i_j\leq q} [W_{i_1}\cap\cdots\cap
W_{i_j}]\sK\\
&=\sum\limits_{S} N_{S}\, [W_S]\sK,
\end{align*}
where each $W_S$ is an intersection $W_{i_1}\cap\cdots\cap
W_{i_j}$, and the integer coefficient $N_S$ accounts for the fact
that $W_S$ can in general be expressed as an intersection of
$W_i$'s in more than one way. Since $W_S$ is itself a coordinate
subspace, $[W_S]\sK$ is easily computed.

Gr\"obner degenerations are used to obtain the double Grothendieck
and double Schubert polynomials for all permutations by Knutson
and Miller \cite{Kn-Mi}, and for vexillary permutations by
Knutson, Miller, and Yong \cite{Kn-Mi-Yo2}, \cite{Kn-Mi-Yo}. In
\cite{Ko-Ra}, \cite{Kr-La}, and \cite{Kr}, only Grassmannians
(i.e., Grassmannian permutations) are studied. However, this
allows the authors to degenerate at the local level, which in this
paper results in the positivity of the restriction formulas.
The methods and results of this paper have been extended to the
case of Symplectic Grassmannians by Kreiman \cite{Kr2} (see also
Ikeda \cite{Ik}), by using a Gr\"obner degeneration of
Ghorpade-Raghavan \cite{Gh-Ra}.

\subsection*{Semistandard Set-Valued Tableaux}

Semistandard set-valued tableaux are generalizations of
semistandard Young tableaux. These objects were introduced by Buch
\cite{Bu}, who used them to give a formula for the linear
structure constants for products of Schubert classes in the
K-theory of the Grassmannian. Buch also expressed Grothendieck
polynomials for Grassmannian permutations in terms of semistandard
set-valued tableaux.  Knutson, Miller, and Yong
\cite{Kn-Mi-Yo2,Kn-Mi-Yo} give several formulas for double
Grothendieck and double Schubert polynomials for vexillary
permutations in terms of flagged set-valued tableaux.

In Sections \ref{s.families_nonint_paths} and
\ref{s.three_equiv_models}, we discuss three equivalent
combinatorial models: certain semistandard Young tableaux,
`families of nonintersecting paths on Young diagrams', and
`subsets of Young diagrams'. Although the three models are equally
suitable for expressing our equivariant cohomology formula, we
find the tableau model to be the simplest one for deriving and
expressing the equivariant K-theory formula. Each $[W_i]\sK$ is
naturally indexed by a semistandard Young tableau $P_i$, and each
$[W_{i_1}\cap\cdots\cap W_{i_j}]\sK$ is naturally indexed by the
`union' $S=P_{i_1}\cup\cdots\cup P_{i_j}$, which is a set-valued
tableau. The $S$ for which $N_S\neq 0$ are precisely those which
are semistandard.

The families of nonintersecting paths which we use appeared first
in Krattenthaler \cite{Kra1,Kra2}
and subsequently in \cite{Ko-Ra},  \cite{Kr-La}, and \cite{Kr}.
The subsets of Young diagrams, which were discovered independently
by Ikeda-Naruse, are similar to RC graphs or reduced pipe dreams
\cite{Be-Bi,Fo-Ki,Kn-Mi} for Grassmannian permutations.

\vspace{1em}

The paper is organized as follows.  In Section \ref{s.summary}, we
state our formulas for $[X_\ga]\sK|_{e_\gb}$ and
$[X_\ga]\sK|_{e_\gb}$. In Section \ref{s.equiv_aff_spaces}, we
present basic definitions and properties of equivariant K-theory
and equivariant cohomology for affine spaces and affine varieties.
In Section \ref{s.class_opp_sch_var}, we give the main arguments
for the proof of our formulas (Proposition
\ref{p.rest_formula_kthry}), omitting the proofs of two lemmas. In
Sections \ref{s.families_nonint_paths} and
\ref{s.three_equiv_models}, we prove the first of these two
lemmas, by translating a result of \cite{Ko-Ra}, \cite{Kr-La}, and
\cite{Kr} into the language of semistandard Young tableaux. In
Section \ref{s.combin_svt}, we prove the second lemma, which
computes $N_S$. \vspace{1em}

\noindent\textbf{Acknowledgements.} I would like to thank W.
Graham, P. Magyar, and M. Shimozono for helpful discussions and
suggestions.

\section{Statement of Results}\label{s.summary}

Let $d$ and $n$ be fixed positive integers, $0<d<n$.  Let
$I_{d,n}$ be the set of $d$-element subsets of $\{1,\ldots,n\}$,
where we always assume the entries of such a subset are listed in
increasing order.  We define the \textbf{complement} of
$\ga=\{\ga(1),\ldots,\ga(d)\}\in I_{d,n}$ by
${\ga'}=\{1,\ldots,n\}\setminus \ga$ and the \textbf{length} of
$\ga$ by $l(\ga)=(\ga(1)-1)+\cdots+(\ga(d)-d)$. Let $J_{d,n}$ be
the set of partitions $\gl=(\gl_1,\ldots,\gl_d)$ such that
$n-d\geq \gl_1\geq\cdots\geq \gl_d\geq 0$. There is a standard
bijection $\gp:I_{d,n}\to J_{d,n}$ given by
$\gp(\{\ga(1),\ldots,\ga(d)\})= (\ga(d)-d\ldots,\ga(1)-1)$. Let
$\ga,\gb\in I_{d,n}$ be fixed.

The \textbf{Grassmannian} $Gr_{d,n}$ is the set of all
$d$-dimensional complex subspaces of $\bC^n$. Let
$\{e_1,\ldots,e_n\}$ be the standard basis for $\bC^n$. Define
$e_\ga=\Span\{e_{\ga(1)},\ldots,e_{\ga(d)}\}\in Gr_{d,n}$.
Consider the opposite standard flag, whose $i$-th space is
$F_i=\Span\{e_n,\ldots,e_{n-i+1}\}$, $i=1,\ldots,n$. The
\textbf{Schubert variety} $X_\ga$  (which is sometimes called an
{\it opposite} Schubert variety) is defined by incidence
relations:
\[
X_\ga=\{V\in Gr_{d,n}\mid \dim (V\cap F_i)\geq \dim (e_{\ga}\cap
F_i),\, i=1,\ldots,n\}.
\]

A \textbf{Young diagram} is a collection of boxes arranged into a
left and top justified array. If the $i$-th row of a diagram has
$\gl_i$ boxes, $i=1,\ldots,r$, then we say that the \textbf{shape}
of the diagram is the partition $\gl=(\gl_1,\ldots,\gl_r)$.  A
\textbf{set-valued tableau} is an assignment of a nonempty set of
positive integers to each box of a diagram. The entries of a
set-valued tableau $S$ are the positive integers in the boxes. If
a positive integer occurs in more than one box of $S$, then we
consider the separate occurrences of the positive integer to be
distinct entries of $S$. A \textbf{Young tableau} is a special
type of set-valued tableau in which each box contains a single
entry.


A set-valued tableau is said to be \textbf{semistandard} if all
entries of any box $B$ are less than or equal to all entries of
the box to the right of $B$ and strictly less than all entries of
the box below $B$.\vspace{1em}

\psset{unit=1.3cm}
\begin{center}
\pspicture(0,0)(5,3)
\psset{linewidth=.02}
\psline(0,3)(5,3) \psline(0,2)(5,2) \psline(0,1)(3,1)
\psline(0,0)(2,0)
\psline(0,0)(0,3) \psline(1,0)(1,3) \psline(2,0)(2,3)
\psline(3,1)(3,3) \psline(4,2)(4,3) \psline(5,2)(5,3)
\rput{0}(.5,2.5){$1$} \rput{0}(1.5,2.5){$2,3$}
\rput{0}(2.5,2.5){$3$} \rput{0}(3.5,2.5){$4,6,7$}
\rput{0}(4.5,2.5){$7,9$}
\rput{0}(0.5,1.5){$2$} \rput{0}(1.5,1.5){$4,5,7$}
\rput{0}(2.5,1.5){$8,9$}
\rput{0}(0.5,0.5){$4,6$}  \rput{0}(1.5,0.5){$8$}
\endpspicture
\end{center}

\begin{center}
Figure 1: A semistandard set-valued tableau
\end{center}
\vspace{1em}

If $\gm=(\gm_1,\ldots,\gm_h)$ is any partition, then a set-valued
tableau $S$ is said to be \textbf{on $\gm$} if, for every entry
$x$ of $S$, $x\leq h$ and
\begin{equation}\label{e.semi_stand_on_mu}
x+c(x)-r(x)\leq \gm_x,
\end{equation}
where $r(x)$ and $c(x)$ are the row and column numbers of the box
containing $x$.  Note that the condition $x\leq h$ is required for
$\gm_x$, and thus (\ref{e.semi_stand_on_mu}), to be well-defined.

\newcommand{\ot}{{1,2}}
\newcommand{\twth}{{2,3}}
\begin{ex} Let $\gl=(2,1)$, $\gm=(4,4,2,1)$. The following
list gives all semistandard set-valued tableaux on $\gm$ of shape
$\gl$: \Yboxdim22pt \Yvcentermath1
\[
\begin{array}{c@{\hspace{1.5em}}c@{\hspace{1.5em}}c@{\hspace{1.5em}}c@{\hspace{1.5em}}c@{\hspace{1.5em}}c}
\young(11,2)&\young(12,2)&\young(11,3)&\young(12,3)&\young(22,3)&\young(\ot 2,3)\\[3em]

\young(1\ot,2)&\young(1\ot,3)&\young(11,\twth)&\young(12,\twth)&\young(1\ot,\twth)&
\end{array}
\]
\end{ex}
\vspace{.5em}

\noindent Denote the set of semistandard set-valued tableaux on
$\gm$ of shape $\gl$ by $\ssvt\lm$ and the set of semistandard
Young tableaux on $\gm$ of shape $\gl$ by $\st\lm$.

\begin{prop}\label{p.rest_formula_kthry}
Let $\gl=\gp(\ga)$, $\gm=\gp(\gb)$. Then

\noindent (i) $\displaystyle
[X_\ga]\sK|_{e_\gb}=(-1)^{l(\ga)}\sum\limits_{S\in\ssvt\lm} \,
\prod\limits_{x\in
S}\left(\frac{t_{\gb(d+1-x)}}{t_{{\gb'}(x+c(x)-r(x))}}-1\right)$.
\vspace{1em}

\noindent (ii) $\displaystyle
[X_\ga]\sH|_{e_\gb}=\sum\limits_{S\in\st\lm}\, \prod\limits_{x\in
S}\left(t_{\gb(d+1-x)}-t_{{\gb'}(x+c(x)-r(x))}\right)$.
\end{prop}

\begin{ex} Consider $Gr_{3,6}$, $\ga=\{1,3,5\}$, $\gb=\{2,5,6\}$.
Then $l(\ga)=3$, $\gp(\ga)=(2,1,0)$, $\gp(\gb)=(3,3,1)$. The
semistandard set-valued tableaux on $\gp(\gb)$ of shape $\gp(\ga)$
are:
\Yboxdim22pt \Yvcentermath1
\[
\young(11,2)\qquad\qquad \young(12,2)\qquad\qquad \young(1\ot,2)
\]
Therefore,
\begin{align*}
[X_\ga]\sK|_{e_\gb}&=-\left[\left(\frac{t_6}{t_1}-1\right)\left(\frac{t_6}{t_3}-1\right)\left(\frac{t_5}{t_1}-1\right)
+\left(\frac{t_6}{t_1}-1\right)\left(\frac{t_4}{t_2}-1\right)\left(\frac{t_5}{t_1}-1\right)\right.\\
&\qquad\qquad +\left.\left(\frac{t_6}{t_1}-1\right)\left(\frac{t_6}{t_3}-1\right)\left(\frac{t_4}{t_2}-1\right)\left(\frac{t_5}{t_1}-1\right)\right].\\
[X_\ga]\sH|_{e_\gb}&=(t_6-t_1)(t_6-t_3)(t_5-t_1)+(t_6-t_1)(t_4-t_2)(t_5-t_1).
\end{align*}

\end{ex}

\begin{rem}\label{r.positivity}
In Proposition \ref{p.rest_formula_kthry}(i), the condition that
$S$ is on $\gm$ implies that each term in the product is of the
form $t_b/t_a-1$, $b>a$, and in Proposition
\ref{p.rest_formula_kthry}(ii), the condition that $S$ is on $\gm$
implies that each term in the product is of the form $t_b-t_a$,
$b>a$. Indeed, for $\gm=\gp(\gb)$, one can show that
$\gm_j=\#\{i\in\{1,\ldots,n-d\}\mid {\gb'}(i)< \gb(d+1-j)\}$,
$j=1,\ldots,d$. Therefore $i\leq\gm_j\iff {\gb'}(i)< \gb(d+1-j)$.
Substituting $i=x+c(x)-r(x)$, $j=x$, we obtain:
(\ref{e.semi_stand_on_mu})$\iff x+c(x)-r(x)\leq \gm_x\iff
{\gb'}(x+c(x)-r(x))<\gb(d+1-x)$.
\end{rem}

\section{Equivariant K-Theory in Affine Spaces}\label{s.equiv_aff_spaces}

The equivariant K-theory $K_T^*(V)$ of an algebraic variety $V$
with a $T$ action is defined to be the Grothendieck group of
equivariant coherent sheaves of $\mathcal{O}_V$ modules. If
$Y\subset V$ is a $T$-stable closed subvariety, then we define
$[Y]\sK$ to be the class of the structure sheaf $\mathcal{O}_Y$ of
$Y$.

In this section we assume that $V$ is the affine space $\bC^m$. In
this case, the notion of coherent sheaves of $\mathcal{O}_V$
modules can be replaced by that of finitely generated $\bC[V]$
modules. If $Y$ is a $T$-stable closed subvariety of $V$, then
$[Y]\sK$ is just $[\bC[Y]]\sK$. We also have that $K_T^*(V)\cong
K_T^*({\bf 0})$, which can be identified with the representation
ring of $T$, $R(T)=\bC[t_1^{\pm 1},\ldots,t_n^{\pm 1}]$.

For any (possibly infinite dimensional) $T$-module $L$, define\\
$\Char(L)\in \bC[[t_1^{\pm 1},\ldots,t_n^{\pm 1}]]$ to be the
character of $L$ under the $T$ action (one also views $\Char(L)$
as the $\mathbb{Z}^n$-graded Hilbert function of $L$, where each
character of $L$ is graded by its $T$-weight). For ${\bf
d}\in\bb{Z}^n$, define $\bC[V](-{\bf d})$ to be $\bC[V]$ with
modified $T$-action: the characters of $\bC[V](-{\bf d})$ are the
same as those of $\bC[V]$, but with weight ${\bf t^d}$ times
greater. We use the standard identification

\begin{equation*}
[\bC[V](-{\bf d})]\sK={\bf t}^{\bf d}=\frac{\Char(\bC[V](-{\bf
d}))}{\Char(\bC[V])}.
\end{equation*}

Let $Y$ be a $T$-stable closed subvariety of $V$.  There is a free
equivariant resolution
\begin{equation*}
0\rightarrow\cE_r\rightarrow\cdots\rightarrow\cE_1\rightarrow\bC[Y]\rightarrow
0,\ \ \hbox{ where }\ \ \cE_i=\bigoplus\limits_{j=1}^{u_i}
\bC[V](-{\bf d}_{ij}).
\end{equation*}
Since
\begin{equation*}
[\cE_i]\sK=\sum\limits_{j=1}^{u_i}[\bC[V](-{\bf
d}_{ij})]\sK=\sum\limits_{j=1}^{u_i}\frac{\Char(\bC[V](-{\bf
d}_{ij}))}{\Char(\bC[V])}=\frac{\Char(\cE_i)}{\Char(\bC[V])},
\end{equation*}
it follows that
\begin{equation}\label{eq.k_poly}
[Y]\sK=\sum\limits_{i=1}^r
(-1)^{i+1}[\cE_i]\sK=\sum\limits_{i=1}^r(-1)^{i+1}\frac{\Char(\cE_i)}{\Char(\bC[V])}
=\frac{\Char(\bC[Y])}{\Char(\bC[V])}.
\end{equation}


\begin{ex}\label{ex.k_poly}
Let $V=\bC^3$, and let $y_1,y_2,y_3\in\bC[V]$ be the standard
coordinate functions on $V$. Suppose that $T=(\bC^*)^4$ acts on
$V$, and hence on $\bC[V]$, and suppose that for ${\bf
t}=\text{diag}(t_1,t_2,t_3,t_4)\in T$,
\begin{equation*}
{\bf t}(y_1)=\frac{t_4}{t_1}\,y_1,\qquad {\bf
t}(y_2)=\frac{t_2}{t_3}\,y_2,\qquad {\bf t}(y_3)=t_1^{-2}\,y_3.
\end{equation*}
Then
\begin{align*}
\Char(\bC[V])&=\sum\limits_{i,j,k=0}^{\infty}
\left(\frac{t_4}{t_1}\right)^i
\left(\frac{t_2}{t_3}\right)^j
\left(t_1^{-2}\right)^k
=\frac{1}{(1-\frac{t_4}{t_1})(1-\frac{t_2}{t_3})(1-t_1^{-2})}.\\
\intertext{Let $Y\subset V$ be the $y_1$-axis. Then}
\Char(\bC[Y])&=\sum\limits_{i=0}^{\infty}
\left(\frac{t_4}{t_1}\right)^i=\frac{1}{(1-\frac{t_4}{t_1})}
=\frac{(1-\frac{t_2}{t_3})(1-t_1^{-2})}{(1-\frac{t_4}{t_1})(1-\frac{t_2}{t_3})(1-t_1^{-2})}.
\end{align*}
Therefore, by (\ref{eq.k_poly}),
$[Y]\sK=(1-\frac{t_2}{t_3})(1-t_1^{-2})$.
\end{ex}

\noindent Let $y_1,\ldots,y_m\in\bC[V]$ be the standard coordinate
functions on $V=\bC^m$. We denote by $V(\{
y_{j_1},\ldots,y_{j_k}\})$ the coordinate subspace of $V$ defined
by the vanishing of $y_{j_1},\ldots,y_{j_k}$.

\begin{lem}\label{l.k_poly_union_coord_subsp}
Let $\chi_i$, $i=1,\ldots,m$ be characters of $T$.  Suppose that
$T$ acts on $V$, and hence on $\bC[V]$, and suppose that ${\bf
t}(y_i)=\chi_i({\bf t}) y_i$, ${\bf t}\in T$,
$i=1,\ldots,m$.\\
\noindent (i)   If $W\subset V$ is the coordinate subspace $W=V(\{
y_{j_1},\ldots,y_{j_k}\})$, then
\begin{equation}\label{eq.k_poly_coord_subsp}
[W]\sK=(1-\chi_{j_1}({\bf t}))\cdots(1-\chi_{j_k}({\bf t})).\\
\end{equation}

\noindent (ii) If $W\subset V$ is the union of coordinate
subspaces $W_1,\ldots,W_q$, then
\begin{equation}\label{e.k_poly_union_coord_subsp} [W]\sK=
\sum\limits_{j=1}^q(-1)^{j+1}\sum\limits_{1\leq i_1<\cdots<i_j\leq
q} [W_{i_1}\cap\cdots\cap W_{i_j}]\sK.
\end{equation}
\end{lem}
\begin{proof}
(i) is an easy generalization of Example \ref{ex.k_poly}.\\

\noindent (ii) By the inclusion-exclusion principle,
\begin{align*}
\Char(W)&=\Char(W_1\cup\cdots\cup W_q)\\
&=\sum\limits_{j=1}^q(-1)^{j+1}\sum\limits_{1\leq
i_1<\cdots<i_j\leq q} \Char(W_{i_1}\cap\cdots\cap W_{i_j})
\end{align*}
The result now follows from (\ref{eq.k_poly}).
\end{proof}
\noindent Note that each $W_{i_1}\cap\cdots\cap W_{i_j}$ in
(\ref{e.k_poly_union_coord_subsp}) is itself a coordinate
subspace, so (\ref{eq.k_poly_coord_subsp}) can be used to compute
its class.
%

\section{The Class of an Opposite Schubert Variety}\label{s.class_opp_sch_var}

The \textbf{Pl\"ucker map} $\pl:Gr_{d,n}\to\mathbb{P}(\wedge^d
\bC^n)$ is defined by $\pl(W)=[w_1\wedge\cdots\wedge w_d]$, where
$\{w_1,\ldots,w_d\}$ is any basis for $W$. It is well known that
$\pl$ is a closed immersion. Thus $Gr_{d,n}$ inherits the
structure of projective variety, as does $X_\ga\subset Gr_{d,n}$.

\subsection*{Reduction to an Affine Variety}

Under the Pl\"ucker map, $e_\gb$ maps to
$[e_{\gb_1}\wedge\cdots\wedge e_{\gb_d}]\in\mathbb{P}(\wedge^d
\bC^n)$. Define $p_\gb$ to be homogeneous (\textbf{Pl\"ucker})
coordinate $[e_{\gb_1}\wedge\cdots\wedge
e_{\gb_d}]^*\in\bC[\mathbb{P}(\wedge^d \bC^n)]$. Let $\cO_\gb$ be
the distinguished open set of $Gr_{d,n}$ defined by $p_\gb\neq 0$.
Then $\cO_\gb$ is isomorphic to the affine space $\bC^{d(n-d)}$,
with $e_\gb$ the origin. Indeed, $\cO_\gb$ can be identified with
the space of matrices in $M_{n\times d}$ in which rows
$\gb_1,\ldots,\gb_d$ are the rows of the $d\times d$ identity
matrix, and rows ${\gb'}_1,\ldots,{\gb'}_{n-d}$ contain arbitrary
elements of $\bC$. Under this identification, the rows of
$\cO_\gb$ are indexed by $\{1,\ldots,n\}$, and the columns by
$\gb$.

\begin{ex}
Let $d=3$, $n=7$, $\gb=\{2,5,7\}$. Then $\gb'=\{1,3,4,6\}$, and
\begin{align*}
\cO_\gb&=\left\{ \left(
\begin{matrix}
y_{12}&y_{15}&y_{17}\\
1&0&0\\
y_{32}&y_{35}&y_{37}\\
y_{42}&y_{45}&y_{47}\\
0&1&0\\
y_{62}&y_{65}&y_{67}\\
0&0&1
\end{matrix}
\right), y_{ab}\in \bC \right\}.
\end{align*}
\end{ex}
\noindent The space $\cO_\gb$ is $T$-stable, and for ${\bf
t}=\hbox{diag}(t_1,\ldots,t_n)\in T$ and coordinate functions
$y_{ab}\in\bC[\cO_\gb]$,
\begin{equation}\label{e.char_O_gb}
{\bf t}(y_{ab})=\frac{t_b}{t_a}\,y_{ab}.
\end{equation}

The equivariant embeddings
$e_\gb\stackrel{j}{\to}\cO_\gb\stackrel{k}{\to} Gr_{d,n}$ induce
homomorphisms
\[
K_T^*(Gr_{d,n})\stackrel{k^*}{\rightarrow}
K_T^*(\cO_\gb)\stackrel{j^*}{\rightarrow} K_T^*(e_\gb).
\]
The map $j^*$ is an isomorphism, identifying $ K_T^*(\cO_\gb)$
with $K_T^*(e_\gb)$. Define $Y_{\ga,\gb}=X_\ga\cap\cO_\gb$. We
have
\begin{equation}\label{e.restrict_to_affine}
[X_\ga]\sK|_{e_\gb}=j^*\circ
k^*([X_\ga]\sK)=j^*([k^{-1}X_\ga]\sK)=j^*([Y\ab]\sK)=[Y\ab]\sK.
\end{equation}
Applying analogous arguments for equivariant cohomology, we obtain
\begin{equation}\label{e.restrict_to_affine_cohom}
[X_\ga]\sH|_{e_\gb}=[Y\ab]\sH.
\end{equation}

\subsection*{Reduction to a Union of Coordinate Subspaces}

Let $\gl=\gp(\ga)$, $\gm=\gp(\gb)$. Let $\svt\lm$ denote the set
of all set-valued tableaux (not necessarily semistandard) of shape
$\gl$ on $\gm$.  For $S\in\svt\lm$, define
\begin{equation}\label{e.W_S_defn}
W_S=V(\{ y_{{\gb'}(x+c(x)-r(x)),\gb(d+1-x)}, x\in S \}),
\end{equation}
a coordinate subspace of $\cO_\gb$. Define
\begin{equation}\label{e.W_ab_defn}
W\ab=\bigcup\limits_{P\in\st\lm}W_P.
\end{equation}
The following lemma, whose proof appears in Section
\ref{s.three_equiv_models}, reduces our problem to computing the
class of a union of coordinate subspaces.
\begin{lem}\label{l.restrict_to_planes}
$[Y\ab]\sK=[W\ab]\sK$
\end{lem}

Let $R$ and $S$ be two set-valued tableaux of shape $\gl$. Define
the \textbf{union} $R\cup S$ to be the set-valued tableau of shape
$\gl$ whose entries in each box are the unions of the entries of
$R$ and the entries of $S$ in that box. If $R$ and $S$ are both on
$\gm$, then $R\cup S$ is on $\gm$, and $W_{R\cup S}=W_R\cap
W_{S}$. We say that $R$ is \textbf{contained in} $S$, and write
$R\subset S$, if each entry in each box of $R$ is also an entry in
the same box of $S$. In this case, if $S$ is on $\gm$, then $R$ is
on $\gm$.

Let $S$ be a semistandard tableau of shape $\gl$. Define $\cP(S)$
to be the set of semistandard Young tableax of shape $\gl$ which
are contained in $S$, and define $q_S=|\cP(S)|$. Define $N_{S,j}$
to be the number of $j$ element subsets of $\cP(S)$ whose unions
equal $S$, and define $N_S=\sum\limits_{j=1}^{q_S}(-1)^{j+1}
N_{S,j}$.

\begin{lem}\label{l.class_union_planes_sch}
$\displaystyle
[W\ab]\sK=\sum\limits_{S\in\svt\lm}N_{S}\prod\limits_{x\in
S}\left(1-\frac{t_{\gb(d+1-x)}}{t_{{\gb'}(x+c(x)-r(x))}}\right)$.
\end{lem}
\begin{proof}
Let $P_1,\ldots,P_q$ be an enumeration of $\st\lm$. Note that for
any $S\in\svt\lm$, $\cP(S)\subset\st\lm$; thus $q_S\leq q$. By
(\ref{e.W_ab_defn}) and Lemma
\ref{l.k_poly_union_coord_subsp}(ii),
\begin{align*}
[W\ab]\sK&=\left[ W_{P_1}\cup\dots\cup W_{P_q}\right]\sK\\
&=\sum\limits_{j=1}^q(-1)^{j+1}\sum\limits_{1\leq
i_1<\cdots<i_j\leq q} [W_{P_{i_1}}\cap\cdots\cap
W_{P_{i_j}}]\sK\\
&=\sum\limits_{j=1}^q(-1)^{j+1}\sum\limits_{1\leq
i_1<\cdots<i_j\leq q} [W_{P_{i_1}\cup\cdots\cup P_{i_j}}]\sK\\
&=\sum\limits_{j=1}^q(-1)^{j+1}\ \sum\limits_{S\in\svt\lm}\
\sum\limits_{P_{i_1}\cup\cdots\cup P_{i_j}=S,\atop 1\leq
i_1<\cdots<i_j\leq q}
  [W_S]\sK\\
&=\sum\limits_{j=1}^q(-1)^{j+1}\ \sum\limits_{S\in\svt\lm}
N_{S,j}\,[W_S]\sK\\
%
&=\sum\limits_{S\in\svt\lm}\ \sum\limits_{j=1}^q(-1)^{j+1}
N_{S,j}\,[W_S]\sK\\
&=\sum\limits_{S\in\svt\lm}\ \sum\limits_{j=1}^{q_S}(-1)^{j+1}
N_{S,j}\,[W_S]\sK\\
&=\sum\limits_{S\in\svt\lm} N_{S}\, [W_S]\sK.
\end{align*}
By (\ref{e.W_S_defn}), (\ref{e.char_O_gb}), and Lemma
\ref{l.k_poly_union_coord_subsp}(i), for each $S\in\svt\lm$,
\begin{equation*}
[W_S]\sK=\prod\limits_{x\in
S}\left(1-\frac{t_{\gb(d+1-x)}}{t_{{\gb'}(x+c(x)-r(x))}}\right).
\end{equation*}
\end{proof}

\noindent For set-valued tableau $S$, define $|S|$ (resp. $\|S\|$)
to be the total number of entries (resp. boxes) of $S$. The proof
of the following Lemma appears in Section \ref{s.combin_svt}.

\begin{lem}\label{l.N_S}
(i) If $S$ is semistandard, then $N_S=(-1)^{|S|+\|S\|}$.\\
(ii) If $S$ is not semistandard, then $N_S=0$.
\end{lem}



\begin{proof}[Proof of Proposition \ref{p.rest_formula_kthry}.]
(i) Combining the preceding reductions and lemmas:
\begin{align*}
[X_\ga]\sK|_{e_\gb}&=[Y\ab]\sK\\
&=[W\ab]\sK\\
&=\sum\limits_{S\in\svt\lm} N_{S}\prod\limits_{x\in
S}\left(1-\frac{t_{\gb(d+1-x)}}{t_{{\gb'}(x+c(x)-r(x))}}\right)\\
&=\sum\limits_{S\in\ssvt\lm} (-1)^{|S|+\|S\|}\prod\limits_{x\in
S}\left(1-\frac{t_{\gb(d+1-x)}}{t_{{\gb'}(x+c(x)-r(x))}}\right)\\
&=(-1)^{l(\ga)}\sum\limits_{S\in\ssvt\lm}
(-1)^{|S|}\prod\limits_{x\in
S}\left(1-\frac{t_{\gb(d+1-x)}}{t_{{\gb'}(x+c(x)-r(x))}}\right)\\
&=(-1)^{l(\ga)}\sum\limits_{S\in\ssvt\lm}\, \prod\limits_{x\in
S}\left(\frac{t_{\gb(d+1-x)}}{t_{{\gb'}(x+c(x)-r(x))}}-1\right).
\end{align*}

\noindent (ii) By (\ref{e.groth_schub}) and \cite{Kn-Mi}, Lemma
1.1.4, $[X\ab]\sH|_{e_\gb}({\bf t})$ equals the sum of the lowest
degree terms of $[X\ab]\sK|_{e_\gb}({\bf 1- t})$. One checks that
the sum of the lowest degree terms of
\begin{equation*}
[X\ab]\sK|_{e_\gb}({\bf 1- t})
=(-1)^{l(\ga)}\sum\limits_{S\in\ssvt\lm}\, \prod\limits_{x\in
S}\left(\frac{1-t_{\gb(d+1-x)}}{1-t_{{\gb'}(x+c(x)-r(x))}}-1\right)
\end{equation*}
equals
\begin{equation*}
\sum\limits_{S\in\st\lm}\, \prod\limits_{x\in
S}\left(t_{\gb(d+1-x)}-t_{{\gb'}(x+c(x)-r(x))}\right).
\end{equation*}
\end{proof}

\section{Families of Nonintersecting Paths on Young
Diagrams}\label{s.families_nonint_paths}

In this section, we introduce a set $\cF\lm$ of families of
nonintersecting paths on Young diagrams, which we then show to be
identical to the set $\cF''\lm$ of families of nonintersecting
paths which appear in \cite{Ko-Ra}, \cite{Kra1}, \cite{Kra2},
\cite{Kr-La}, and \cite{Kr}, despite the fact that $\cF\lm$ and
$\cF''\lm$ are defined quite differently. Thus one can view the
result of this section as giving an alternate way of defining or
expressing the path families of \cite{Ko-Ra}, \cite{Kra1},
\cite{Kra2}, \cite{Kr-La}, and \cite{Kr}.  In terms of this
alternate definition, one can more easily see the equivalence
between the path families and the other combinatorial models in
Section \ref{s.three_equiv_models}.

Let $\gm$ be a partition and $D_\gm$ the corresponding Young
diagram.   We denote by $(i,j)$ the box of $D_\gm$ at row $i$
(from the top), column $j$ (from the left).  We impose the
following order on boxes of $D_\gm$: if $(i,j),(k,l)\in D_\gm$,
then $(i,j)\leq (k,l)$ if $i\leq k$ and $j\leq l$.

A \textbf{path on} $D_\gm$ is a path of contiguous boxes of the
Young diagram $D_\gm$ which
\begin{itemize}
\item[(i)] moves only up or to the right, and

\item[(ii)] begins on the lowest box of a column and ends on the
rightmost box of a row.
\end{itemize}
Note that a path on $D_\gm$ may consist of only one box. In this
case, the box must be a lower right corner of $D_\gm$.  We define
the \textbf{greatest lower bound} of a path $P$, or $\glb(P)$, to
be the greatest lower bound of all the boxes of $P$. Explicitly,
if the left endpoint of $P$ is the box $(i,j)$ and the right
endpoint is the box $(k,l)$, then $\glb(P)$ is the box $(k,j)\in
D_\gm$. In particular, the endpoints of $P$ determine $\glb(P)$.
We impose the following order on paths on $D_\gm$: if $P,P'$ are
paths on $D_\gm$, then $P\leq P'$ if $\glb(P)\leq\glb(P')$.

Denote by $\cF_\gm$ the set of families of nonintersecting paths
on $D_\gm$. Let $F\in\cF_\gm$. Define the \textbf{support} of $F$,
$\supp(F)$, to be the set of all boxes in all paths of $F$.

\begin{ex}\label{e.path_family} A family $F=\{P_1,P_2,P_3,P_4\}$ of nonintersecting paths on
$D_\gm$, $\gm=(9,9,9,9,8,8,6,6,3,1)$. We have $\glb(P_1)=(1,1)$,
$\glb(P_2)=(3,5)$, $\glb(P_3)=(5,7)$, $\glb(P_4)=(9,2)$, and
$P_1\leq P_2\leq P_3$, $P_1\leq P_4$.
\psset{unit=.35cm}
\begin{center}
\setlength{\arraycolsep}{.6cm} \pspicture(0,0)(9,11)
\psset{linewidth=.015}
\psline(0,9)(9,9) \psline(0,10)(9,10)
\psline(0,8)(9,8) \psline(0,7)(9,7) \psline(0,6)(9,6)
\psline(0,5)(8,5) \psline(0,4)(8,4) \psline(0,3)(6,3)
\psline(0,2)(6,2) \psline(0,1)(3,1) \psline(0,0)(1,0)
\psline(0,0)(0,10) \psline(1,0)(1,10) \psline(2,1)(2,10)
\psline(3,1)(3,10) \psline(4,2)(4,10) \psline(5,2)(5,10)
\psline(6,2)(6,10) \psline(7,4)(7,10) \psline(8,4)(8,10)
\psline(9,6)(9,10)
\rput(-.5,.5){\tiny $10$}
\rput(-.5,1.5){\tiny $9$}
\rput(-.5,2.5){\tiny $8$}
\rput(-.5,3.5){\tiny $7$}
\rput(-.5,4.5){\tiny $6$}
\rput(-.5,5.5){\tiny $5$}
\rput(-.5,6.5){\tiny $4$}
\rput(-.5,7.5){\tiny $3$}
\rput(-.5,8.5){\tiny $2$}
\rput(-.5,9.5){\tiny $1$}
\rput(.5,10.5){\tiny $1$}
\rput(1.5,10.5){\tiny $2$}
\rput(2.5,10.5){\tiny $3$}
\rput(3.5,10.5){\tiny $4$}
\rput(4.5,10.5){\tiny $5$}
\rput(5.5,10.5){\tiny $6$}
\rput(6.5,10.5){\tiny $7$}
\rput(7.5,10.5){\tiny $8$}
\rput(8.5,10.5){\tiny $9$}
\psset{linewidth=.3}
\psline(.5,.35)(.5,4.5)(2.5,4.5)(2.5,5.5)(3.5,5.5)(3.5,7.5)(4.5,7.5)(4.5,8.5)(8.5,8.5)(8.5,9.65)
\psline(4.5,2.35)(4.5,3.5)(5.5,3.5)(5.5,6.5)(6.5,6.5)(6.5,7.5)(8.65,7.5)
\psline(1.35,1.5)(2.65,1.5)
\psline(6.35,4.5)(7.5,4.5)(7.5,5.65)
\rput(2.6,6.5){\tiny $P_1$}
\rput(4.6,4.5){\tiny $P_2$}
\rput(6.6,5.4){\tiny $P_3$}
\rput(2.5,2.4){\tiny $P_4$}
\endpspicture

\end{center}
\end{ex}

\begin{lem}\label{l.supp_dets_family}$\supp(F)$ uniquely determines the  paths of $F$.
\end{lem}
\begin{proof}
In other words, if $\supp(F')=\supp(F)$, then $F'=F$. Our proof is
by decreasing induction on the number $m$ of paths of $F$. If $P$
is a minimal path of $F$, then one sees that $P$ must also be a
path of $F'$. In particular, the result for $m=1$ holds.  If
$m>1$, then $\supp(F'\setminus
\{P\})=\supp(F')\setminus\supp(\{P\})=\supp(F)\setminus\supp(\{P\})=\supp(F\setminus
\{P\})$, and $F\setminus \{P\}$ has $m-1$ paths.  Thus by
induction $F'\setminus \{P\}=F\setminus \{P\}$, which completes
the proof.
\end{proof}

Let $\gl\leq\gm$, i.e., $\gl_i\leq\gm_i$
for all $i$. Then $D_\gl$ naturally embeds in $D_\gm$ in such
a way that both Young diagrams share the same top left corner.  In
this way, $D_\gl$ can be viewed as a subset of $D_\gm$.

\eject 

\begin{lem}\label{l.flmtop_exists}
There exists $F\in\cF_\gm$ such that $\supp(F)=D_\gm\setminus
D_\gl$.
\end{lem}
\begin{proof}
We use decreasing induction on the number of boxes of
$D_\gm\setminus D_\gl$. One can form a path $P$ on $D_\gm$
consisting of boxes $(i,j)\in D_\gm\setminus D_\gl$ such that
$(i-1,j-1)\not\in D_\gm\setminus D_\gl$.  This path moves along an
upper left boundary of $D_\gm\setminus D_\gl$.  Now $D'=D_\gl\cup
P$ is a Young diagram contained in $D_\gm$, and $D_\gm\setminus
D'$ has fewer boxes than $D_\gm\setminus D_\gl$. Thus by induction
there exists $F'\in\cF_\gm$ such that $\supp(F')=D_\gm\setminus
D'=(D\gm\setminus D_\gl)\setminus P$. Set $F=F'\cup\{P\}$.
\end{proof}
\noindent By Lemma \ref{l.supp_dets_family}, the family $F$ of
Lemma \ref{l.flmtop_exists} is uniquely determined. We call this
family the \textbf{top family on} $\mathbf{D_\gm}$ \textbf{for}
$\gl$ and denote it by $F\lm\tp$.

\begin{ex} The family $F\lm\tp$, where $\gm=(7,6,6,6,3,3)$,
$\gl=(3,2)$. \vspace{5pt} \psset{unit=.35cm}
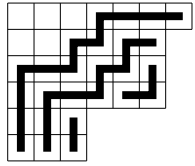
\begin{center}
\setlength{\arraycolsep}{.6cm} \pspicture(0,0)(7,6)
\psset{linewidth=.015}
\psline(0,6)(7,6) \psline(0,5)(7,5) \psline(0,4)(6,4)
\psline(0,3)(6,3) \psline(0,2)(6,2) \psline(0,1)(3,1)
\psline(0,0)(3,0)
\psline(0,0)(0,6) \psline(1,0)(1,6) \psline(2,0)(2,6)
\psline(3,0)(3,6) \psline(4,2)(4,6) \psline(5,2)(5,6)
\psline(6,2)(6,6) \psline(7,5)(7,6)

\psset{linewidth=.3}
\psline(.5,.35)(.5,3.5)(2.5,3.5)(2.5,4.5)(3.5,4.5)(3.5,5.5)(6.65,5.5)
\psline(1.5,.35)(1.5,2.5)(3.5,2.5)(3.5,3.5)(4.5,3.5)(4.5,4.5)(5.65,4.5)
\psline(2.5,.35)(2.5,1.65) \psline(4.35,2.5)(5.5,2.5)(5.5,3.65)
\endpspicture
\end{center}
\end{ex}
\vspace{5pt}

If $F\in F_\gm$, then we define $\tw(F)=\{\glb(P)\mid P\text{ is a
path of }F\}$, and we call this the \textbf{twist} of $F$. (In
Example \ref{e.path_family},
$\tw(F)=\{(1,1),(3,5),(5,7),(9,2)\}$.)
\begin{lem}\label{l.twist_unique}
If $\tw(F_{\gn,\gm}\tp)=\tw(F\lm\tp)$, then $\gn=\gl$.
\end{lem}
\begin{proof}
By decreasing induction on $\gm$. Let $b\in D_{\gm}$ be a maximal
element of $\tw(F\lm\tp)$, and thus also of $\tw(F_{\gn,\gm}\tp)$.
Then $b=\glb(P)$ for some maximal path $P$ of $F\lm\tp$, which
runs along the lower-right boundary of $D_\gm$. Likewise,
$b=\glb(P')$ for some maximal path $P'$ of $F_{\gn,\gm}\tp$, which
runs along the lower-right boundary of $D_\gm$. Since $P'$ and $P$
are paths with the same greatest lower bound, and they both run
along the lower right boundary of $D_\gm$, $P'=P$.

Now $D_\gm\setminus P=D_{\gm'}$, for some $\gm'<\gm$. We have that
$F_{\gl,\gm'}\tp=F\lm\tp\setminus\{P\}$,
$F_{\gn,\gm'}\tp=F\lm\tp\setminus\{P\}$, and thus
$\tw(F_{\gn,\gm'})=\tw(F_{\gn,\gm})\setminus\{b\}
=\tw(F_{\gl,\gm})\setminus\{b\}=\tw(F_{\gl,\gm'})$. Therefore
$\gn=\gl$ follows by induction.
\end{proof}

Suppose that $D_\gm$ contains the four boxes $(i,j)$, $(i+1,j)$,
$(i,j+1)$, and $(i+1,j+1)$ (which make up a square), and that $F$
contains $(i+1,j)$, $(i,j+1)$, and $(i+1,j+1)$, but not $(i,j)$.
One checks that these three boxes must lie on the same path $P$ of
$F$. We apply a \textbf{ladder move} to $F$ by altering $P$ as
follows: the box $(i+1,j+1)$ of $P$ is removed and replaced with
the box $(i,j)$, thus obtaining a new path $P'$ on $D_\gm$. The
resulting path $P'$, combined with the paths of $F$ other than
$P$, form a new family of nonintersecting paths $F'$ on $D_\gm$.
We denote this ladder move by $F\to F'$. Note that a ladder move
is invertible. We call the inverse of a ladder move a
\textbf{reverse ladder move}.

\eject

\begin{ex}
Let $F$ be the following family of nonintersecting paths on
$D_\gm$, $\gm=(6,5,5,4,4,1)$: \vspace{5pt} \psset{unit=.35cm}
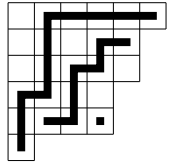
\begin{center}
\setlength{\arraycolsep}{.6cm} \pspicture(0,0)(6,6)
\psset{linewidth=.015}
\psline(0,6)(6,6) \psline(0,5)(6,5) \psline(0,4)(5,4)
\psline(0,3)(5,3) \psline(0,2)(4,2) \psline(0,1)(4,1)
\psline(0,0)(1,0)
\psline(0,0)(0,6) \psline(1,0)(1,6) \psline(2,1)(2,6)
\psline(3,1)(3,6) \psline(4,1)(4,6) \psline(5,3)(5,6)
\psline(6,5)(6,6)

\psset{linewidth=.3}
\psline(.5,.35)(.5,2.5)(1.5,2.5)(1.5,5.5)(5.65,5.5)
\psline(1.35,1.5)(2.5,1.5)(2.5,3.5)(3.5,3.5)(3.5,4.5)(4.65,4.5)
\psline(3.35,1.5)(3.65,1.5)
\endpspicture
\end{center}

\noindent The following two ladder moves can be applied to $F$:
\vspace{5pt} \psset{unit=.35cm}
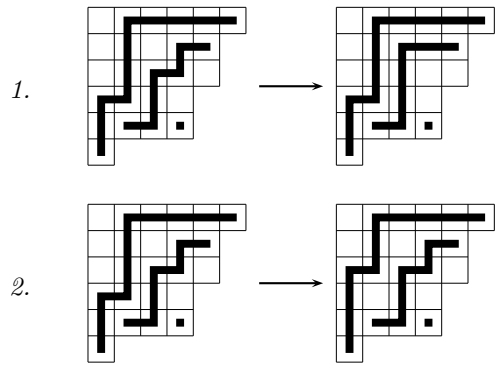
\begin{center}
1. $ \setlength{\arraycolsep}{.6cm}
\begin{array}{cc}
\rnode{a} { \pspicture(0,0)(6,6)
\psset{linewidth=.015}
\psline(0,6)(6,6) \psline(0,5)(6,5) \psline(0,4)(5,4)
\psline(0,3)(5,3) \psline(0,2)(4,2) \psline(0,1)(4,1)
\psline(0,0)(1,0)
\psline(0,0)(0,6) \psline(1,0)(1,6) \psline(2,1)(2,6)
\psline(3,1)(3,6) \psline(4,1)(4,6) \psline(5,3)(5,6)
\psline(6,5)(6,6)

\psset{linewidth=.3}
\psline(.5,.35)(.5,2.5)(1.5,2.5)(1.5,5.5)(5.65,5.5)
\psline(1.35,1.5)(2.5,1.5)(2.5,3.5)(3.5,3.5)(3.5,4.5)(4.65,4.5)
\psline(3.35,1.5)(3.65,1.5)
\endpspicture
}
&\rnode{b} {\pspicture(0,0)(6,6)
\psset{linewidth=.015}
\psline(0,6)(6,6) \psline(0,5)(6,5) \psline(0,4)(5,4)
\psline(0,3)(5,3) \psline(0,2)(4,2) \psline(0,1)(4,1)
\psline(0,0)(1,0)
\psline(0,0)(0,6) \psline(1,0)(1,6) \psline(2,1)(2,6)
\psline(3,1)(3,6) \psline(4,1)(4,6) \psline(5,3)(5,6)
\psline(6,5)(6,6)

\psset{linewidth=.3}
\psline(.5,.35)(.5,2.5)(1.5,2.5)(1.5,5.5)(5.65,5.5)
\psline(1.35,1.5)(2.5,1.5)(2.5,4.5)(4.65,4.5)
\psline(3.35,1.5)(3.65,1.5)
\endpspicture
}
\end{array}
$ \ncline[nodesep=5pt]{->}{a}{b} \vspace{1em}

2. $\setlength{\arraycolsep}{.6cm}
\begin{array}{cc}
\rnode{a} { \pspicture(0,0)(6,6)
\psset{linewidth=.015}
\psline(0,6)(6,6) \psline(0,5)(6,5) \psline(0,4)(5,4)
\psline(0,3)(5,3) \psline(0,2)(4,2) \psline(0,1)(4,1)
\psline(0,0)(1,0)
\psline(0,0)(0,6) \psline(1,0)(1,6) \psline(2,1)(2,6)
\psline(3,1)(3,6) \psline(4,1)(4,6) \psline(5,3)(5,6)
\psline(6,5)(6,6)

\psset{linewidth=.3}
\psline(.5,.35)(.5,2.5)(1.5,2.5)(1.5,5.5)(5.65,5.5)
\psline(1.35,1.5)(2.5,1.5)(2.5,3.5)(3.5,3.5)(3.5,4.5)(4.65,4.5)
\psline(3.35,1.5)(3.65,1.5)
\endpspicture
}
&\rnode{b} {\pspicture(0,0)(6,6)
\psset{linewidth=.015}
\psline(0,6)(6,6) \psline(0,5)(6,5) \psline(0,4)(5,4)
\psline(0,3)(5,3) \psline(0,2)(4,2) \psline(0,1)(4,1)
\psline(0,0)(1,0)
\psline(0,0)(0,6) \psline(1,0)(1,6) \psline(2,1)(2,6)
\psline(3,1)(3,6) \psline(4,1)(4,6) \psline(5,3)(5,6)
\psline(6,5)(6,6)

\psset{linewidth=.3}
\psline(.5,.35)(.5,3.5)(1.5,3.5)(1.5,5.5)(5.65,5.5)
\psline(1.35,1.5)(2.5,1.5)(2.5,3.5)(3.5,3.5)(3.5,4.5)(4.65,4.5)
\psline(3.35,1.5)(3.65,1.5)
\endpspicture
}
\end{array}
$ \ncline[nodesep=5pt]{->}{a}{b}
\end{center}
\end{ex}
\vspace{5pt}

A \textbf{family of nonintersecting paths on} $\mathbf{D_\gm}$
\textbf{for} $\gl$ is an element of $\cF_\gm$ which can be
obtained by applying a succession of ladder moves to $F\lm\tp$.
The set of all families of nonintersecting paths on
$\mathbf{D_\gm}$ for $\gl$ is denoted by $\cF\lm$.

\begin{ex} The following diagram shows $\cF\lm$, as well as all
possible ladder moves, where $\gm=(4,4,3,3,1)$, $\gl=(2,1)$.
\vspace{5pt} \psset{unit=.35cm}
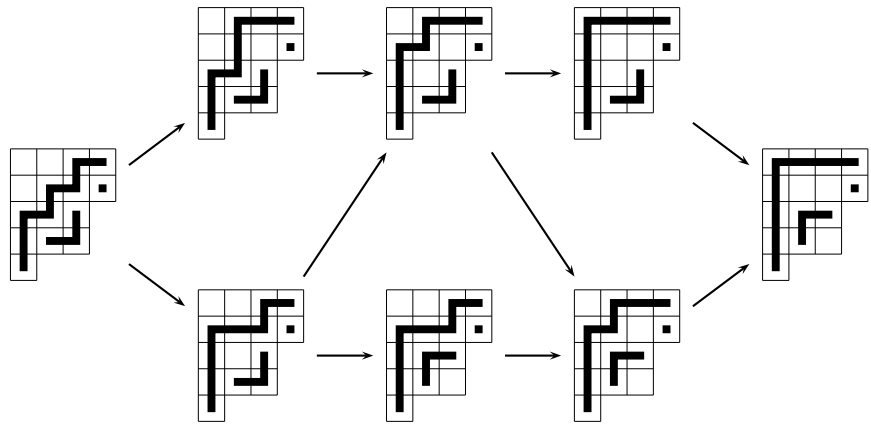
\begin{center}
$ \setlength{\arraycolsep}{.55cm}
\begin{array}{lcccc}
&\rnode{b}{ \pspicture(0,0)(4,5)
\psset{linewidth=.015}
\psline(0,5)(4,5) \psline(0,4)(4,4) \psline(0,3)(4,3)
\psline(0,2)(3,2) \psline(0,1)(3,1) \psline(0,0)(1,0)
\psline(0,0)(0,5) \psline(1,0)(1,5) \psline(2,1)(2,5)
\psline(3,1)(3,5) \psline(4,3)(4,5)
\psset{linewidth=.3}
\psline(.5,.35)(.5,2.5)(1.5,2.5)(1.5,4.5)(3.65,4.5)
\psline(1.35,1.5)(2.5,1.5)(2.5,2.65) \psline(3.35,3.5)(3.65,3.5)
\endpspicture
}
&\rnode{c}{ \pspicture(0,0)(4,5)
\psset{linewidth=.015}
\psline(0,5)(4,5) \psline(0,4)(4,4) \psline(0,3)(4,3)
\psline(0,2)(3,2) \psline(0,1)(3,1) \psline(0,0)(1,0)
\psline(0,0)(0,5) \psline(1,0)(1,5) \psline(2,1)(2,5)
\psline(3,1)(3,5) \psline(4,3)(4,5)
\psset{linewidth=.3}
\psline(.5,.35)(.5,3.5)(1.5,3.5)(1.5,4.5)(3.65,4.5)
\psline(1.35,1.5)(2.5,1.5)(2.5,2.65) \psline(3.35,3.5)(3.65,3.5)
\endpspicture
}
&\rnode{d}{ \pspicture(0,0)(4,5)
\psset{linewidth=.015}
\psline(0,5)(4,5) \psline(0,4)(4,4) \psline(0,3)(4,3)
\psline(0,2)(3,2) \psline(0,1)(3,1) \psline(0,0)(1,0)
\psline(0,0)(0,5) \psline(1,0)(1,5) \psline(2,1)(2,5)
\psline(3,1)(3,5) \psline(4,3)(4,5)
\psset{linewidth=.3} \psline(.5,.35)(.5,4.5)(3.65,4.5)
\psline(1.35,1.5)(2.5,1.5)(2.5,2.65) \psline(3.35,3.5)(3.65,3.5)
\endpspicture
} &
\\
\rnode{a}{ \pspicture(0,0)(4,5)
\psset{linewidth=.015}
\psline(0,5)(4,5) \psline(0,4)(4,4) \psline(0,3)(4,3)
\psline(0,2)(3,2) \psline(0,1)(3,1) \psline(0,0)(1,0)
\psline(0,0)(0,5) \psline(1,0)(1,5) \psline(2,1)(2,5)
\psline(3,1)(3,5) \psline(4,3)(4,5)
\psset{linewidth=.3}
\psline(.5,.35)(.5,2.5)(1.5,2.5)(1.5,3.5)(2.5,3.5)(2.5,4.5)(3.65,4.5)
\psline(1.35,1.5)(2.5,1.5)(2.5,2.65) \psline(3.35,3.5)(3.65,3.5)
\endpspicture
} &&& &\rnode{e}{ \pspicture(0,0)(4,5)
\psset{linewidth=.015}
\psline(0,5)(4,5) \psline(0,4)(4,4) \psline(0,3)(4,3)
\psline(0,2)(3,2) \psline(0,1)(3,1) \psline(0,0)(1,0)
\psline(0,0)(0,5) \psline(1,0)(1,5) \psline(2,1)(2,5)
\psline(3,1)(3,5) \psline(4,3)(4,5)
\psset{linewidth=.3} \psline(.5,.35)(.5,4.5)(3.65,4.5)
\psline(1.5,1.35)(1.5,2.5)(2.65,2.5) \psline(3.35,3.5)(3.65,3.5)
\endpspicture
}
\\
&\rnode{B}{ \pspicture(0,0)(4,5)
\psset{linewidth=.015}
\psline(0,5)(4,5) \psline(0,4)(4,4) \psline(0,3)(4,3)
\psline(0,2)(3,2) \psline(0,1)(3,1) \psline(0,0)(1,0)
\psline(0,0)(0,5) \psline(1,0)(1,5) \psline(2,1)(2,5)
\psline(3,1)(3,5) \psline(4,3)(4,5)
\psset{linewidth=.3}
\psline(.5,.35)(.5,3.5)(2.5,3.5)(2.5,4.5)(3.65,4.5)
\psline(1.35,1.5)(2.5,1.5)(2.5,2.65) \psline(3.35,3.5)(3.65,3.5)
\endpspicture
}
&\rnode{C}{ \pspicture(0,0)(4,5)
\psset{linewidth=.015}
\psline(0,5)(4,5) \psline(0,4)(4,4) \psline(0,3)(4,3)
\psline(0,2)(3,2) \psline(0,1)(3,1) \psline(0,0)(1,0)
\psline(0,0)(0,5) \psline(1,0)(1,5) \psline(2,1)(2,5)
\psline(3,1)(3,5) \psline(4,3)(4,5)
\psset{linewidth=.3}
\psline(.5,.35)(.5,3.5)(2.5,3.5)(2.5,4.5)(3.65,4.5)
\psline(1.5,1.35)(1.5,2.5)(2.65,2.5) \psline(3.35,3.5)(3.65,3.5)
\endpspicture
}
&\rnode{D}{ \pspicture(0,0)(4,5)
\psset{linewidth=.015}
\psline(0,5)(4,5) \psline(0,4)(4,4) \psline(0,3)(4,3)
\psline(0,2)(3,2) \psline(0,1)(3,1) \psline(0,0)(1,0)
\psline(0,0)(0,5) \psline(1,0)(1,5) \psline(2,1)(2,5)
\psline(3,1)(3,5) \psline(4,3)(4,5)
\psset{linewidth=.3}
\psline(.5,.35)(.5,3.5)(1.5,3.5)(1.5,4.5)(3.65,4.5)
\psline(1.5,1.35)(1.5,2.5)(2.65,2.5) \psline(3.35,3.5)(3.65,3.5)
\endpspicture
} &
\ncline[nodesep=5pt]{->}{a}{b} \ncline[nodesep=5pt]{->}{b}{c}
\ncline[nodesep=5pt]{->}{c}{d} \ncline[nodesep=5pt]{->}{d}{e}
\ncline[nodesep=5pt]{->}{a}{B} \ncline[nodesep=5pt]{->}{B}{C}
\ncline[nodesep=5pt]{->}{C}{D} \ncline[nodesep=5pt]{->}{D}{e}
\ncline[nodesep=5pt]{->}{B}{c} \ncline[nodesep=5pt]{->}{c}{D}
\end{array}
$
\end{center}
\end{ex}
\vspace{5pt}

\subsection*{The Set $\cF''\lm$ of Path Families}

We say that a set of boxes $S\subset D_\gm$ is a \textbf{twisted
chain} if for any two boxes $p$, $p'$ in $S$:
\begin{itemize}
\item[1.] $p$ and $p'$ lie on different rows and different columns
of $D_\gm$.
\item[2.] Either $p<p'$, $p'<p$, or $\{p,p'\}$ has no upper bound
in $D_\gm$.
\end{itemize}

\begin{ex}  The shaded boxes form a twisted chain in $D_\gm$.
\psset{unit=.35cm}
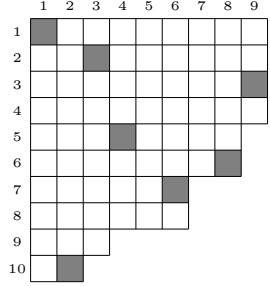
\begin{center}
\setlength{\arraycolsep}{.6cm} \pspicture(0,0)(9,11)
\psset{linewidth=.015}
\psframe*[linecolor=gray](0,9)(1,10)
\psframe*[linecolor=gray](2,8)(3,9)
\psframe*[linecolor=gray](8,7)(9,8)
\psframe*[linecolor=gray](3,5)(4,6)
\psframe*[linecolor=gray](5,3)(6,4)
\psframe*[linecolor=gray](7,4)(8,5)
\psframe*[linecolor=gray](1,0)(2,1)
\psline(0,9)(9,9) \psline(0,10)(9,10) \psline(0,8)(9,8)
\psline(0,7)(9,7) \psline(0,6)(9,6) \psline(0,5)(8,5)
\psline(0,4)(8,4) \psline(0,3)(6,3) \psline(0,2)(6,2)
\psline(0,1)(3,1) \psline(0,0)(2,0)
\psline(0,0)(0,10) \psline(1,0)(1,10) \psline(2,0)(2,10)
\psline(3,1)(3,10) \psline(4,2)(4,10) \psline(5,2)(5,10)
\psline(6,2)(6,10) \psline(7,4)(7,10) \psline(8,4)(8,10)
\psline(9,6)(9,10)
\rput(-.5,.5){\tiny $10$} \rput(-.5,1.5){\tiny $9$}
\rput(-.5,2.5){\tiny $8$} \rput(-.5,3.5){\tiny $7$}
\rput(-.5,4.5){\tiny $6$} \rput(-.5,5.5){\tiny $5$}
\rput(-.5,6.5){\tiny $4$} \rput(-.5,7.5){\tiny $3$}
\rput(-.5,8.5){\tiny $2$} \rput(-.5,9.5){\tiny $1$}
\rput(.5,10.5){\tiny $1$} \rput(1.5,10.5){\tiny $2$}
\rput(2.5,10.5){\tiny $3$} \rput(3.5,10.5){\tiny $4$}
\rput(4.5,10.5){\tiny $5$} \rput(5.5,10.5){\tiny $6$}
\rput(6.5,10.5){\tiny $7$} \rput(7.5,10.5){\tiny $8$}
\rput(8.5,10.5){\tiny $9$}
\endpspicture
\end{center}
Note that, for example, $\{(10,2),(5,4)\}$ has least upper bound
$(10,4)$, which is not in $D_\gm$.
\end{ex}
\vspace{5pt}

It follows from the definitions that if $F\in \cF_\gm$, then
$\tw(F)$ is a twisted chain of $D_\gm$. Define
$$\cF'\lm=\{F\in\cF\gm\mid\tw(F)=\tw(F\lm\tp)\}.$$ Untimately we
will be interested in the set $\cF''\lm$ of families of
nonintersecting paths on $D_\gm$, which is defined below. We
introduce $\cF'\lm$ in order to break up the proof that
$\cF\lm=\cF''\lm$ into two parts: $\cF\lm=\cF'\lm$ and
$\cF'\lm=\cF''\lm$.
\begin{lem}\label{l.path1_path2}
$\cF'\lm=\cF\lm$.
\end{lem}
\begin{proof}
We remark that applying a ladder move to any family in $\cF\lm$
does not alter its twist. Thus ladder moves preserve $\cF'\lm$.
Also $F\lm\tp\in\cF'\lm$.  Since $\cF\lm$ is by definition the
smallest subset of $\cF_\gm$ preserved by ladder moves which
contains $F\lm\tp$, $\cF\lm\subset\cF'\lm$.

Next suppose that $F\in\cF'\lm$. By applying a succession of
reverse ladder moves to $F$, $F_{\gn,\gm}\tp$ can be obtained, for
some $\gn\leq\gm$. Since reverse latter moves do not alter twist,
$\tw(F_{\gn,\gm}\tp)=\tw(F)=\tw(F_{\gl,\gm}\tp)$. By Lemma
\ref{l.twist_unique}, $\gn=\gl$. Thus $F$ can be obtained from
$F\lm\tp$ by applying ladder moves, which implies $F\in\cF\lm$.
Hence $\cF'\lm\subset\cF\lm$.
\end{proof}

Let $\gl\leq\gm$, and define $\ga=\gp^{-1}(\gl)$,
$\gb=\gp^{-1}(\gm)$. In \cite{Kr}, it is shown that there exists a
unique twisted chain $S\lm=\{(x_1,y_1),\ldots,(x_t,y_t)\}\subset
D_\gm$ such that
$\ga=\gb\setminus\{\gb(d+1-x_1),\ldots,\gb(d+1-x_t)\}\cup\{\gb'(y_1),\ldots,\gb'(y_t)\}$.
Define $$\cF''\lm=\{F\in\cF_\gm\mid\tw(F)=S\lm\}.$$ It follows
from the definitions that $\cF''\lm=\cF'\lm$ if and only if
$\tw(F\lm\tp)=S\lm$. We first prove this for the case where $\gl$
and $\gm$ are such that $D_\gm\setminus D_\gl$ consists of a
single path.

\begin{lem}\label{l.part_dtuple_1}
Let $P$ be a path which moves along the lower right border of
$D_\gm$, so that $D_\gm\setminus \{P\}=D_{\gn}$, for some
$\gn<\gm$. Then $\tw(F_{\gn,\gm}\tp)=S_{\gn,\gm}$.
\end{lem}

\begin{ex} Let $d=5$, $n=9$.  Let $\gm=(4,4,3,3,1)$, and let $P$ be the path on
$D_\gm$ shown below.  Then $D_\gm\setminus P=D_\gn$, where
$\gn=(4,2,2,1,1)$.  We have $F_{\gn,\gm}\tp=\{P\}$, and
$\tw(F_{\gn,\gm}\tp)=\glb(P)=(2,2)$. Let
$\gb=\gp^{-1}(\gm)=\{2,5,6,8,9\}$,
$\ga=\gp^{-1}(\gn)=\{2,3,5,6,9\}$. Then $\ga=\gb\setminus 8\cup
3=\gb\setminus \gb(5+1-2)\cup\gb'(2)$, so $S_{\gn,\gm}=(2,2)$,
which agrees with $\tw(F_{\gn,\gm}\tp)$.

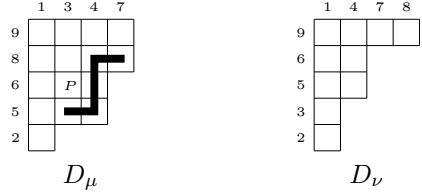
\psset{unit=.35cm}
\setlength{\arraycolsep}{.6cm}
\begin{equation*}
\begin{array}{ccc}
\pspicture(0,0)(4,6) \psset{linewidth=.015}
\psline(0,5)(4,5) \psline(0,4)(4,4) \psline(0,3)(4,3)
\psline(0,2)(3,2) \psline(0,1)(3,1) \psline(0,0)(1,0)
\psline(0,0)(0,5) \psline(1,0)(1,5) \psline(2,1)(2,5)
\psline(3,1)(3,5) \psline(4,3)(4,5)
\rput(-.5,.5){\tiny $2$} \rput(-.5,1.5){\tiny $5$}
\rput(-.5,2.5){\tiny $6$} \rput(-.5,3.5){\tiny $8$}
\rput(-.5,4.5){\tiny $9$}
\rput(.5,5.5){\tiny $1$} \rput(1.5,5.5){\tiny $3$}
\rput(2.5,5.5){\tiny $4$} \rput(3.5,5.5){\tiny $7$}
\psset{linewidth=.3}
\psline(1.35,1.5)(2.5,1.5)(2.5,3.5)(3.65,3.5)
\rput(1.6,2.5){\tiny $P$}
\endpspicture
&&
\pspicture(0,0)(4,6) \psset{linewidth=.015}
\psline(0,5)(4,5) \psline(0,4)(4,4) \psline(0,3)(2,3)
\psline(0,2)(2,2) \psline(0,1)(1,1) \psline(0,0)(1,0)
\psline(0,0)(0,5) \psline(1,0)(1,5) \psline(2,2)(2,5)
\psline(3,4)(3,5) \psline(4,4)(4,5)
\rput(-.5,.5){\tiny $2$} \rput(-.5,1.5){\tiny $3$}
\rput(-.5,2.5){\tiny $5$} \rput(-.5,3.5){\tiny $6$}
\rput(-.5,4.5){\tiny $9$}
\rput(.5,5.5){\tiny $1$} \rput(1.5,5.5){\tiny $4$}
\rput(2.5,5.5){\tiny $7$} \rput(3.5,5.5){\tiny $8$}
\endpspicture\\
D_\gm && D_\gn
\end{array}
\end{equation*}
\end{ex}

\begin{proof}[Proof of Lemma \ref{l.part_dtuple_1}]
For $i=1,\ldots,d$, define $\ol{i}=d+1-i$; thus if a box is on row
$i$ of $D_\gm$ counting from the top, then it is on row $\ol{i}$
counting from the bottom. Let $(x,y)=\glb(P)=\tw(F_{\gn,\gm}\tp)$.
Let $\gb=\gp^{-1}(\gm)$ and $\ga=\gb\setminus\gb(\ol{x})\cup
\gb'(y)$. By definition, $S_{\gp(\ga),\gm}=(x,y)$. We shall obtain
the result by showing that $\gp(\ga)=\gn$.

Note that
\begin{equation*}
D_\gm=\{(u,v)\in\{1,\ldots,d\}\times\{1,\ldots,n-d\}\mid
\gb(\ol{u})>\gb'(v)\}.
\end{equation*}
Let $w=\min\{u\mid \gb(u)>\gb'(y)\}$, and note that $w\leq
\ol{x}$. We have
$$\ga=\{\gb(1),\dots,\gb(w-1),\gb'(y),\gb(w),\dots,
\gb(\ol{x}-1),\widehat{\gb(\ol{x})},\gb(\ol{x}+1),\dots,\gb(d)\}.$$
Therefore,
\begin{equation*}
\ga(i)=
\begin{cases}
\gb(i),& i<w\ \text{ or }\ i>\ol{x}\\
\gb(i-1), & w<i\leq \ol{x}\\
\gb'(y), & i=w.
\end{cases}
\end{equation*}
We have
\begin{align*}
(\gp(\ga))_i&=\ga(\ol{i})-(\ol{i})\\
& =
\begin{cases}
\gb(\ol{i})-\ol{i},& \ol{i}<w\ \text{ or }\ \ol{i}>\ol{x}\\
\gb(\ol{i}-1)-\ol{i}, & w<\ol{i}\leq \ol{x}\\
\gb'(y)-\ol{w}, & \ol{i}=w
\end{cases}\\
&=
\begin{cases}
\gm_i,& \ol{i}<w\ \text{ or }\ \ol{i}>\ol{x}\\
\gm_{i+1}-1, & w<\ol{i}\leq \ol{x}-1\\
y-1, & \ol{i}=w
\end{cases}\\
&=\gn_i,
\end{align*}
where we use the facts that $\gb(\ol{i})-(\ol{i})=\gm_i$,
$\gb(\ol{i}-1)-(\ol{i})=\gb(\ol{i+1})-(\ol{i+1})-1=\gm_{i+1}-1$,
and $\gb'(y)-y=\ol{w}-1$.
%
\end{proof}

\begin{lem}\label{l.part_dtuple_2}
For any $\gl\leq\gm$, $\tw(F\lm\tp)=S\lm$.
\end{lem}
\begin{proof}
Define $\gb=\gp^{-1}(\gm)$, $\ga=\gp^{-1}(\gl)$, and let
$\{(x_1,y_1),\ldots,(x_t,y_t)\}=\tw(F\lm\tp)$. We wish to show
that $\ga=\gb\setminus\{\gb(d+1-x_1),\ldots,\gb(d+1-x_t)\}
\cup\{\gb'(y_1),\ldots,\gb'(y_t)\}$.

We use decreasing induction on $t$, the number of paths in
$F\lm\tp$. For $t=1$, the result is identically Lemma
\ref{l.part_dtuple_1}. Suppose that $t>1$.  Let $P$ be a maximal
path of $F\lm\tp$.  Then $P$ runs along the lower right border of
$D_\gm$, so $D_\gm\setminus P=D_\gn$, for some $\gn<\gm$. Define
$\gc=\gp^{-1}(\gn)$, and assume that $(x_t,y_t)=\glb(P)$
(re-indexing if necessary). By Lemma \ref{l.part_dtuple_1},
\begin{equation}\label{e.tt1}
\gc=\gb\setminus\gb(d+1-x_t)\cup\gb'(y_t).
\end{equation}
Now $F_{\gl,\gn}\tp=F\lm\tp\setminus\{P\}$, whose twist is equal
to $\{(x_1,y_1),\ldots,(x_{t-1},y_{t-1})\}$. Thus by induction,
\begin{equation}\label{e.tt2}
\ga=\gc\setminus\{\gb(d+1-x_1),\ldots,\gb(d+1-x_{t-1})\}
\cup\{\gb'(y_1),\ldots,\gb'(y_{t-1})\}.
\end{equation}
Combining (\ref{e.tt1}) and (\ref{e.tt2}), we arrive at the
result.
\end{proof}

\noindent From Lemma \ref{l.part_dtuple_2}, we obtain
\begin{lem}\label{l.path2_path3}
$\cF''\lm=\cF'\lm$.
\end{lem}
\noindent Now Lemmas \ref{l.path1_path2} and \ref{l.path2_path3}
imply
\begin{lem}\label{l.path1_path3}
$\cF''\lm=\cF\lm$.
\end{lem}

\section{From Path Families to Semistandard Young Tableaux}\label{s.three_equiv_models}
\Yboxdim11pt

Lemma \ref{l.restrict_to_planes} follows from a result of
\cite{Ko-Ra}, \cite{Kr-La}, and \cite{Kr} which gives an
equivariant Gr\"obner degeneration of a Schubert variety in the
neighborhood of a $T$-fixed point to a union of coordinate
subspaces. However, whereas the result of \cite{Ko-Ra},
\cite{Kr-La}, and \cite{Kr} is expressed in terms of the set
$\cF''\lm$ of families of nonintersecting paths on Young diagrams,
our results, and in particular Lemma \ref{l.restrict_to_planes},
require the semistandard Young tableaux $\st\lm$.  Most the
previous section and this one are taken up in showing the
equivalence between these two combinatorial models. In the
previous section we showed that $\cF''\lm=\cF\lm$.  In this
section, we introduce a new model, `subsets of Young diagrams'. We
will be interested in a certain set of subsets of the Young
diagram $D_\gm$, which we denote by $\cD\lm$. We show that
$\cF\lm$ and $\st\lm$ are both equivalent to $\cD\lm$, and thus
are equivalent to eachother. We summarize the steps we take in
showing the equivalence between $\cF''\lm$ and $\st\lm$ as
follows:
\begin{equation*}
\cF''\lm=\cF\lm\longleftrightarrow\cD\lm\longleftrightarrow\st\lm
\end{equation*}
The subsets $\cD\lm$, which were discovered independently by
Ikeda-Naruse, are similar to RC graphs or reduced pipe dreams
\cite{Be-Bi,Fo-Ki,Kn-Mi} for Grassmannian permutations.

\subsection*{Subsets of Young Diagrams}

Let $\gm$ be a partition and $D_\gm$ the corresponding Young
diagram. Lemma \ref{l.supp_dets_family} tells us that a family of
nonintersecting paths on $D_\gm$ is completely characterized by
its support. This suggests that in order to study families of
nonintersecting paths on $D_\gm$, it suffices to study the
supports of these path families (or the complements in $D_\gm$ of
their supports). This motivates the following definitions.

A \textbf{subset of} $\mathbf{D_\gm}$ is a set of boxes in
$D_\gm$. Let $D$ be a subset of $D_\gm$. Suppose that $D_\gm$
contains the four boxes $(i,j)$, $(i+1,j)$, $(i,j+1)$, and
$(i+1,j+1)$ (which make up a square), but of these four boxes, $D$
only contains $(i,j)$. Then a \textbf{ladder move} removes $(i,j)$
from $D$ and replaces it with $(i+1,j+1)$, thereby obtaining a new
subset $D'$ of $D_\gm$. We denote this ladder move by $D\to D'$.
Note that a ladder move is invertible.  We call the inverse of a
ladder move a \textbf{reverse ladder move}.

\begin{ex} Let $D$ be the following subset of $D_\gm$,
$\gm=(6,5,5,4,4,1)$ (boxes of $D$ are shaded): \vspace{5pt}
\psset{unit=.35cm}
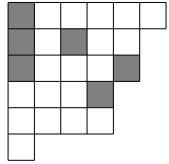
\begin{center}
\setlength{\arraycolsep}{.6cm} \pspicture(0,0)(6,6)
\psframe*[linecolor=gray](0,5)(1,6)
\psframe*[linecolor=gray](0,4)(1,5)
\psframe*[linecolor=gray](0,3)(1,4)
\psframe*[linecolor=gray](2,4)(3,5)
\psframe*[linecolor=gray](3,2)(4,3)
\psframe*[linecolor=gray](4,3)(5,4) \psset{linewidth=.03}
\psline(0,6)(6,6) \psline(0,5)(6,5) \psline(0,4)(5,4)
\psline(0,3)(5,3) \psline(0,2)(4,2) \psline(0,1)(4,1)
\psline(0,0)(1,0)
\psline(0,0)(0,6) \psline(1,0)(1,6) \psline(2,1)(2,6)
\psline(3,1)(3,6) \psline(4,1)(4,6) \psline(5,3)(5,6)
\psline(6,5)(6,6)
\endpspicture
\end{center}

\noindent The following two ladder moves can be applied to $F$:
\vspace{5pt}
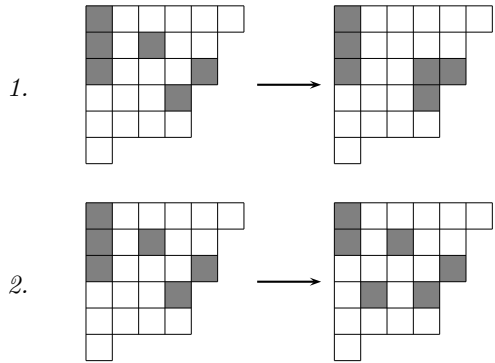
\begin{center}
\psset{unit=.35cm} \setlength{\arraycolsep}{.6cm} 1.
$\begin{array}{cc} \rnode{a}{ \pspicture(0,0)(6,6)
\psframe*[linecolor=gray](0,5)(1,6)
\psframe*[linecolor=gray](0,4)(1,5)
\psframe*[linecolor=gray](0,3)(1,4)
\psframe*[linecolor=gray](2,4)(3,5)
\psframe*[linecolor=gray](3,2)(4,3)
\psframe*[linecolor=gray](4,3)(5,4) \psset{linewidth=.03}
\psline(0,6)(6,6) \psline(0,5)(6,5) \psline(0,4)(5,4)
\psline(0,3)(5,3) \psline(0,2)(4,2) \psline(0,1)(4,1)
\psline(0,0)(1,0)
\psline(0,0)(0,6) \psline(1,0)(1,6) \psline(2,1)(2,6)
\psline(3,1)(3,6) \psline(4,1)(4,6) \psline(5,3)(5,6)
\psline(6,5)(6,6)
\endpspicture
}
&\rnode{b}{ \pspicture(0,0)(6,6)
\psframe*[linecolor=gray](0,5)(1,6)
\psframe*[linecolor=gray](0,4)(1,5)
\psframe*[linecolor=gray](0,3)(1,4)
\psframe*[linecolor=gray](3,3)(4,4)
\psframe*[linecolor=gray](3,2)(4,3)
\psframe*[linecolor=gray](4,3)(5,4) \psset{linewidth=.03}
\psline(0,6)(6,6) \psline(0,5)(6,5) \psline(0,4)(5,4)
\psline(0,3)(5,3) \psline(0,2)(4,2) \psline(0,1)(4,1)
\psline(0,0)(1,0)
\psline(0,0)(0,6) \psline(1,0)(1,6) \psline(2,1)(2,6)
\psline(3,1)(3,6) \psline(4,1)(4,6) \psline(5,3)(5,6)
\psline(6,5)(6,6)
\endpspicture
}
\end{array}
$ \ncline[nodesep=5pt]{->}{a}{b} \vspace{1em}

2. $\setlength{\arraycolsep}{.6cm}\begin{array}{cc} \rnode{a}{
\pspicture(0,0)(6,6)
\psframe*[linecolor=gray](0,5)(1,6)
\psframe*[linecolor=gray](0,4)(1,5)
\psframe*[linecolor=gray](0,3)(1,4)
\psframe*[linecolor=gray](2,4)(3,5)
\psframe*[linecolor=gray](3,2)(4,3)
\psframe*[linecolor=gray](4,3)(5,4) \psset{linewidth=.03}
\psline(0,6)(6,6) \psline(0,5)(6,5) \psline(0,4)(5,4)
\psline(0,3)(5,3) \psline(0,2)(4,2) \psline(0,1)(4,1)
\psline(0,0)(1,0)
\psline(0,0)(0,6) \psline(1,0)(1,6) \psline(2,1)(2,6)
\psline(3,1)(3,6) \psline(4,1)(4,6) \psline(5,3)(5,6)
\psline(6,5)(6,6)
\endpspicture
}
&\rnode{b}{ \pspicture(0,0)(6,6)
\psframe*[linecolor=gray](0,5)(1,6)
\psframe*[linecolor=gray](0,4)(1,5)
\psframe*[linecolor=gray](1,2)(2,3)
\psframe*[linecolor=gray](2,4)(3,5)
\psframe*[linecolor=gray](3,2)(4,3)
\psframe*[linecolor=gray](4,3)(5,4) \psset{linewidth=.03}
\psline(0,6)(6,6) \psline(0,5)(6,5) \psline(0,4)(5,4)
\psline(0,3)(5,3) \psline(0,2)(4,2) \psline(0,1)(4,1)
\psline(0,0)(1,0)
\psline(0,0)(0,6) \psline(1,0)(1,6) \psline(2,1)(2,6)
\psline(3,1)(3,6) \psline(4,1)(4,6) \psline(5,3)(5,6)
\psline(6,5)(6,6)
\endpspicture
}
\end{array}
$ \ncline[nodesep=5pt]{->}{a}{b}
\end{center}
\end{ex}
\vspace{5pt}

Let $\gl$ be a partition with $\gl\leq\gm$, i.e., $\gl_i\leq\gm_i$
for all $i$. Embed the Young diagram $D_\gl$ in $D_\gm$ in such a
way that both $D_\gl$ and $D_\gm$ share the same top left corners.
The subset of $D_\gm$ consisting of all boxes in this embedded
Young diagram is called the \textbf{top subset of}
$\mathbf{D_\gm}$ \textbf{for} $\gl$ and denoted by $D\lm\tp$.

\begin{ex} The subset $D\lm\tp$, where
$\gm=(7,6,6,6,3,3)$, $\gl=(3,2)$. \vspace{3pt} \psset{unit=.35cm}
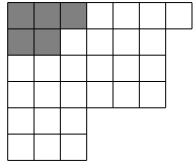
\begin{center}
\setlength{\arraycolsep}{.6cm} \pspicture(0,0)(7,6)
\psframe*[linecolor=gray](0,5)(1,6)
\psframe*[linecolor=gray](0,4)(1,5)
\psframe*[linecolor=gray](1,5)(2,6)
\psframe*[linecolor=gray](1,4)(2,5)
\psframe*[linecolor=gray](2,5)(3,6) \psset{linewidth=.03}
\psline(0,6)(7,6) \psline(0,5)(7,5) \psline(0,4)(6,4)
\psline(0,3)(6,3) \psline(0,2)(6,2) \psline(0,1)(3,1)
\psline(0,0)(3,0)
\psline(0,0)(0,6) \psline(1,0)(1,6) \psline(2,0)(2,6)
\psline(3,0)(3,6) \psline(4,2)(4,6) \psline(5,2)(5,6)
\psline(6,2)(6,6) \psline(7,5)(7,6)
\endpspicture
\end{center}
\end{ex}
\vspace{5pt}

\noindent A \textbf{subset of} $\mathbf{D_\gm}$ \textbf{for} $\gl$
is a subset of $D_\gm$ which can be obtained by applying a
succession of ladder moves to the top subset of $D_\gm$ for $\gl$.
The set of all subsets of $D_\gm$ for $\gl$ is denoted by
$\cD\lm$.

\begin{ex}\label{ex.subsets_of_D_gm} The following diagram shows $\cD\lm$,
as well as all possible ladder moves, where $\gm=(4,4,3,3,1)$,
$\gl=(2,1)$. \vspace{2pt} \psset{unit=.35cm}
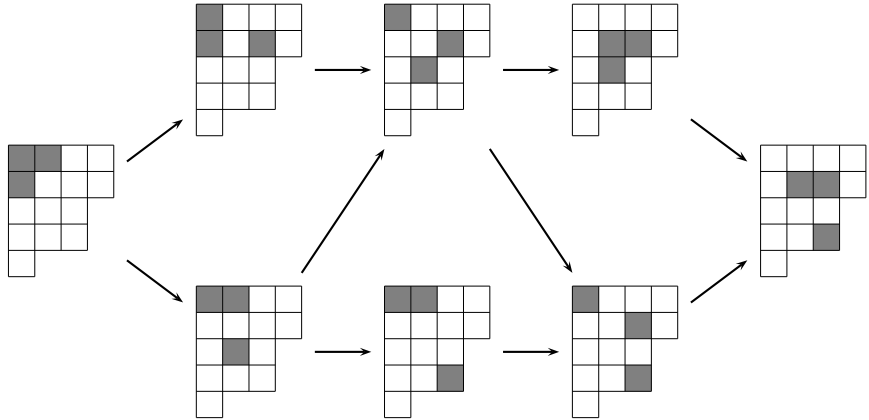
\begin{center}
$ \setlength{\arraycolsep}{.55cm}
\begin{array}{lcccc}
&\rnode{b}{ \pspicture(0,0)(4,5)
\psframe*[linecolor=gray](0,4)(1,5)
\psframe*[linecolor=gray](2,3)(3,4)
\psframe*[linecolor=gray](0,3)(1,4) \psset{linewidth=.03}
\psline(0,5)(4,5) \psline(0,4)(4,4) \psline(0,3)(4,3)
\psline(0,2)(3,2) \psline(0,1)(3,1) \psline(0,0)(1,0)
\psline(0,0)(0,5) \psline(1,0)(1,5) \psline(2,1)(2,5)
\psline(3,1)(3,5) \psline(4,3)(4,5)
\endpspicture
}
&\rnode{c}{ \pspicture(0,0)(4,5)
\psframe*[linecolor=gray](0,4)(1,5)
\psframe*[linecolor=gray](2,3)(3,4)
\psframe*[linecolor=gray](1,2)(2,3) \psset{linewidth=.03}
\psline(0,5)(4,5) \psline(0,4)(4,4) \psline(0,3)(4,3)
\psline(0,2)(3,2) \psline(0,1)(3,1) \psline(0,0)(1,0)
\psline(0,0)(0,5) \psline(1,0)(1,5) \psline(2,1)(2,5)
\psline(3,1)(3,5) \psline(4,3)(4,5)
\endpspicture
}
&\rnode{d}{ \pspicture(0,0)(4,5)
\psframe*[linecolor=gray](1,3)(2,4)
\psframe*[linecolor=gray](2,3)(3,4)
\psframe*[linecolor=gray](1,2)(2,3) \psset{linewidth=.03}
\psline(0,5)(4,5) \psline(0,4)(4,4) \psline(0,3)(4,3)
\psline(0,2)(3,2) \psline(0,1)(3,1) \psline(0,0)(1,0)
\psline(0,0)(0,5) \psline(1,0)(1,5) \psline(2,1)(2,5)
\psline(3,1)(3,5) \psline(4,3)(4,5)
\endpspicture
} &
\\
\rnode{a}{ \pspicture(0,0)(4,5)
\psframe*[linecolor=gray](0,4)(1,5)
\psframe*[linecolor=gray](1,4)(2,5)
\psframe*[linecolor=gray](0,3)(1,4) \psset{linewidth=.03}
\psline(0,5)(4,5) \psline(0,4)(4,4) \psline(0,3)(4,3)
\psline(0,2)(3,2) \psline(0,1)(3,1) \psline(0,0)(1,0)
\psline(0,0)(0,5) \psline(1,0)(1,5) \psline(2,1)(2,5)
\psline(3,1)(3,5) \psline(4,3)(4,5)
\endpspicture
} &&& &\rnode{e}{ \pspicture(0,0)(4,5)
\psframe*[linecolor=gray](1,3)(2,4)
\psframe*[linecolor=gray](2,3)(3,4)
\psframe*[linecolor=gray](2,1)(3,2) \psset{linewidth=.03}
\psline(0,5)(4,5) \psline(0,4)(4,4) \psline(0,3)(4,3)
\psline(0,2)(3,2) \psline(0,1)(3,1) \psline(0,0)(1,0)
\psline(0,0)(0,5) \psline(1,0)(1,5) \psline(2,1)(2,5)
\psline(3,1)(3,5) \psline(4,3)(4,5)
\endpspicture
}
\\
&\rnode{B}{ \pspicture(0,0)(4,5)
\psframe*[linecolor=gray](0,4)(1,5)
\psframe*[linecolor=gray](1,4)(2,5)
\psframe*[linecolor=gray](1,2)(2,3) \psset{linewidth=.03}
\psline(0,5)(4,5) \psline(0,4)(4,4) \psline(0,3)(4,3)
\psline(0,2)(3,2) \psline(0,1)(3,1) \psline(0,0)(1,0)
\psline(0,0)(0,5) \psline(1,0)(1,5) \psline(2,1)(2,5)
\psline(3,1)(3,5) \psline(4,3)(4,5)
\endpspicture
}
&\rnode{C}{ \pspicture(0,0)(4,5)
\psframe*[linecolor=gray](0,4)(1,5)
\psframe*[linecolor=gray](1,4)(2,5)
\psframe*[linecolor=gray](2,1)(3,2) \psset{linewidth=.03}
\psline(0,5)(4,5) \psline(0,4)(4,4) \psline(0,3)(4,3)
\psline(0,2)(3,2) \psline(0,1)(3,1) \psline(0,0)(1,0)
\psline(0,0)(0,5) \psline(1,0)(1,5) \psline(2,1)(2,5)
\psline(3,1)(3,5) \psline(4,3)(4,5)
\endpspicture
}
&\rnode{D}{ \pspicture(0,0)(4,5)
\psframe*[linecolor=gray](0,4)(1,5)
\psframe*[linecolor=gray](2,3)(3,4)
\psframe*[linecolor=gray](2,1)(3,2) \psset{linewidth=.03}
\psline(0,5)(4,5) \psline(0,4)(4,4) \psline(0,3)(4,3)
\psline(0,2)(3,2) \psline(0,1)(3,1) \psline(0,0)(1,0)
\psline(0,0)(0,5) \psline(1,0)(1,5) \psline(2,1)(2,5)
\psline(3,1)(3,5) \psline(4,3)(4,5)
\endpspicture
} &
\ncline[nodesep=5pt]{->}{a}{b} \ncline[nodesep=5pt]{->}{b}{c}
\ncline[nodesep=5pt]{->}{c}{d} \ncline[nodesep=5pt]{->}{d}{e}
\ncline[nodesep=5pt]{->}{a}{B} \ncline[nodesep=5pt]{->}{B}{C}
\ncline[nodesep=5pt]{->}{C}{D} \ncline[nodesep=5pt]{->}{D}{e}
\ncline[nodesep=5pt]{->}{B}{c} \ncline[nodesep=5pt]{->}{c}{D}
\end{array}
$
\end{center}
\end{ex}

\subsection*{Semistandard Young Tableaux}
\Yboxdim13.2pt

A \textbf{ladder move} on a semistandard Young tableau is an
operation which increments one of the entries of the tableau by 1
and results in a semistandard Young tableau.  Note that a ladder
move is invertible. We call the inverse of a ladder move a
\textbf{reverse ladder move}. Recall that a semistandard Young
tableau is \textbf{on} $\gm$ if each of its entries satisfies
(\ref{e.semi_stand_on_mu}). Let $\gl\leq\gm$. Define the
\textbf{top semistandard Young tableau on $\gm$ of shape $\gl$} to
be the (semistandard) Young tableau of shape $\gl$ whose $i$-th
row is filled with $i$'s. This definition does not depend on
$\gm$.

The set of all semistandard Young tableaux on $\gm$ of shape
$\gl$, $\st\lm$, is precisely the set of semistandard Young
tableaux on $\gm$ which can be obtained by applying sequences of
ladder moves to the top semistandard Young tableaux on $\gm$ of
shape $\gl$. This follows from the facts that (i) by applying a
sequence of reverse ladder moves to any semistandard Young tableau
of shape $\gl$, the top semistandard Young tableaux on $\gm$ of
shape $\gl$ can be obtained, and (ii) reverse ladder moves
preserve (\ref{e.semi_stand_on_mu}).

\begin{ex}  The following diagram shows $\st\lm$, as well as all possible ladder moves,
where $\gm=(4,4,3,3,1)$, $\gl=(2,1)$.
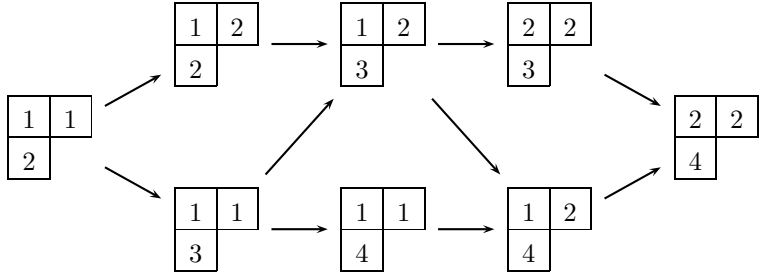
\vspace{2pt}
\begin{center}
\Yboxdim16pt
$ \setlength{\arraycolsep}{.55cm}
\begin{array}{ccccc}
&\rnode{b}{\young(12,2)} &\rnode{c}{\young(12,3)}
&\rnode{d}{\young(22,3)} &
\\
\rnode{a}{\young(11,2)} &&& &\rnode{e}{\young(22,4)}
\\
&\rnode{B}{\young(11,3)} &\rnode{C}{\young(11,4)}
&\rnode{D}{\young(12,4)} &

\ncline[nodesep=5pt]{->}{a}{b}
\ncline[nodesep=5pt]{->}{b}{c}
\ncline[nodesep=5pt]{->}{c}{d}
\ncline[nodesep=5pt]{->}{d}{e}
\ncline[nodesep=5pt]{->}{a}{B}
\ncline[nodesep=5pt]{->}{B}{C}
\ncline[nodesep=5pt]{->}{C}{D}
\ncline[nodesep=5pt]{->}{D}{e}
\ncline[nodesep=5pt]{->}{B}{c}
\ncline[nodesep=5pt]{->}{c}{D}
\end{array}
$
\end{center}
\end{ex}
\vspace{5pt}

\subsection*{The Equivalences $\cF\lm\longleftrightarrow\cD\lm$ and $\cD\lm\longleftrightarrow\st\lm$}


Recall that
\begin{align*}
\cF_\gm&=\hbox{ the set of families of nonintersecting paths on
}D_\gm\\
\cF\lm&=\hbox{ the set of families of nonintersecting paths on
}D_\gm \hbox{ for }\gl\\
F\lm\tp &=\hbox{ the top family of nonintersecting paths on }D_\gm
\hbox{ for }\gl\\[1em]
\cD_\gm&=\hbox{ the set of subsets of }D_\gm\\
\cD\lm&=\hbox{ the set of subsets of }D_\gm\hbox{ for }\gl\\
D\lm\tp &=\hbox{ the top subset of }D_\gm\hbox{ for }\gl\\[1em]
\st\lm&=\hbox{ the set of semistandard Young tableaux on
}\gm\hbox{ of shape
}\gl\\
P\lm\tp&=\hbox{ the top semistandard Young tableaux on }\gm\hbox{
of shape }\gl
\end{align*}
Define $h:\cF\lm\to\cD\lm$ by $h(F)=D_\gm\setminus \supp(F)$. Then
$h$ maps $F\lm\tp$ to $D\lm\tp$ and commutes with ladder moves.
Therefore it restricts to a map from $\cF\lm$ to $\cD\lm$.
Injectivity of $h$ follows immediately from Lemma
\ref{l.supp_dets_family}. To show surjectivity, assume that
$h(F)=D$, for some $D\in \cD_\gm$, $F\in \cF_\gm$. Suppose that a
ladder move is applied to $D$ to obtain $D'$.  Then there is a
corresponding ladder move which when applied to $F$ yields $F'$
such that $h(F')=D'$. Thus surjectivity of $h:\cF\lm\to \cD\lm$
follows by induction on number of ladder moves.

\vspace{1em}


We next give a bijection $g$ from $\cD\lm$ to $\st\lm$. Let $D\in
\cD\lm$. By applying a sequence of reverse ladder moves to $D$,
$D\lm\tp$ can be obtained.  This sequence of reverse ladder moves
takes the box $(x,y)$ of $D$ to some box $(i_{x,y},j_{x,y})$ of
$D\lm\tp$. Note that $(i_{x,y},j_{x,y})$ depends only on $x$ and
$y$ (and $D$), and not on the sequence of reverse ladder moves.
Let $g(D)$ be the tableau of shape $\gl$ which, for each box
$(x,y)$ of $D$, contains entry $x$ in box $(i_{x,y},j_{x,y})$. The
definitions of ladder and reverse ladder moves on subsets of
$D_\gm$ imply that $g(D)$ is semistandard. Also, since $(x,y)\in
D_\gm$, $y\leq \gm(x)$; thus $y-x=j_{x,y}-i_{x,y}$ implies
$x+j_{x,y}-i_{x,y}\leq \gm(x)$. Therefore $g(D)$ is on $\gm$.

To show that $g$ is bijective, we give the inverse map  $f$ from
$\st\lm$ to $\cD\lm$.  Let $P\in\st\lm$, and let $f(P)$ be the
subset of $D_{\gm}$ which, for each entry $x\in P$, contains box
$(x,x-r(x)+c(x))$ of $D_\gm$ (where recall that $(r(x),c(x))$ is
the box of $x$ in $P$). Since $P$ is on $\gm$, $(x,x-r(x)+c(x))\in
D_\gm$; thus $f(P)$ is indeed in $\cD_\gm$. Let $P\lm\tp$ denote
the top semistandard Young tableaux on $\gm$ of shape $\gl$.
Because $f(P\lm\tp)=D\lm\tp$ and $f$ commutes with ladder moves,
$f(P)\in\cD\lm$.

\subsection*{Proof of Lemma \ref{l.restrict_to_planes}}

Let $\gl=\gp(\ga)$, $\gm=\gp(\gb)$. Let $\{v_{a,b}\mid
(a,b)\in{\gb'}\times\gb\}\subset \cO_\gb$ denote the basis dual to
the basis of linear forms $\{y_{a,b}\mid
(a,b)\in{\gb'}\times\gb\}\subset \cO_\gb^*$.  For $F\in\cF_{\gm}$,
define
\begin{equation*}
W_F=\Span(\{v_{{\gb'}(z),\gb(d+1-x)}\mid (x,z)\in \supp(F)\}\cup
\{v_{a,b}\mid (a,b)\in{\gb'}\times\gb, a>b\}).
\end{equation*}
In \cite{Ko-Ra,Kr-La,Kr}, an explicit equivariant bijection (which
is called the bounded RSK in \cite{Kr}) is constructed from
$\bC[\bigcup_{F\in\cF''\lm}W_F]$ to $\bC[Y\ab]$. Thus in light of
Lemma \ref{l.path1_path3},
\begin{equation}\label{e.ml1}
\Char(\bC[Y\ab])=\Char\left(\bC\left[\bigcup_{F\in\cF''\lm}W_F\right]\right)
=\Char\left(\bC\left[\bigcup_{F\in\cF\lm}W_F\right]\right).
\end{equation}

For $x\in\{1,\ldots,d\}$, $z\in\{1,\ldots,n-d\}$, we have that
$(x,z)\in D_\gm\iff z\leq\gm_x\iff {\gb'}(z)<\gb(d+1-x)$ (see
Remark \ref{r.positivity}).  Thus $\{(a,b)\in{\gb'}\times\gb\mid
a<b\}$ can be expressed as $\{({\gb'}(z),\gb(d+1-x))\mid (x,z)\in
D_\gm\}$. Let $F\in\cF\lm$, $D=h(F)$, and $P=g(D)$. Since
$\supp(F)$ and $D$ are complements in $D_\gm$,
\begin{align*} {\gb'}\times\gb
&=\{({\gb'}(z),\gb(d+1-x))\mid (x,z)\in \supp(F)\}\,\dot{\cup}\,
\{(a,b)\in{\gb'}\times\gb\mid a>b\}\\
&\qquad \dot{\cup}\,\{({\gb'}(z),\gb(d+1-x))\mid (x,z)\in D\}.
\end{align*}
Therefore
\begin{align*}
W_F&=V(\{y_{{\gb'}(z),\gb(d+1-x)}\mid (x,z)\in D\})\\
&=V(\{y_{{\gb'}(x-r(x)+c(x)),\gb(d+1-x)}\mid x\in P\})\\
&=W_P.
\end{align*}
Consequently,
$\bigcup_{F\in\cF\lm}W_F=\bigcup_{P\in\st\lm}W_P=W\ab$, and thus
\begin{equation}\label{e.ml2}
\Char\left(\bC\left[\bigcup_{F\in\cF\lm}W_F\right]\right)=
\Char\left(\bC\left[\bigcup_{P\in\st\lm}W_P\right]\right)=\Char(W\ab).
\end{equation}
Combining (\ref{e.ml1}) and (\ref{e.ml2}), we obtain
$\Char(\bC[Y\ab])=\Char(\bC[W\ab])$. By (\ref{eq.k_poly}),
$[Y\ab]\sK=[W\ab]\sK$.


\section{Computing $N_S$}\label{s.combin_svt}

In this section we prove Lemma \ref{l.N_S}.

For a set-valued tableau $S$, we denote by $S_{i,j}$ the set of
entries of $S$ which are contained in the box with row and column
numbers $i$ and $j$ respectively, and we denote by
$S_{i,j,1},\ldots,S_{i,j,r}$ the entries of $S_{i,j}$, which we
assume are listed in increasing order. We define $\cM_k(S)$ to be
the $k$-element subsets of $\sstab(S)$, and $\cN_k(S)$ the
$k$-element subsets of $\sstab(S)$ whose unions equals $S$. By
definition, $N_{S,k}=|\cN_k(S)|$.

Define a \textbf{generalized set-valued tableau} to be the
assignment of a {\it possibly empty} set of positive integers to
each box of a Young diagram. If $S$ and $R$ are two set-valued
tableaux of the same shape, then the difference $S\setminus R$ is
defined to be the generalized set-valued tableau of the same shape
as $S$ (or $R$) with $(S\setminus R)_{i,j}=S_{i,j}\setminus
R_{i,j}$ for all $i,j$. Recall that we write $R\subset S$ to
indicate that $R_{i,j}\subset S_{i,j}$ for all $i,j$.  If
$R\subset S$, and $S$ is clear from the context, then we also
denote $S\setminus R$ by $\ol{R}$. The only generalized set-valued
tableaux which appear in this section in which boxes may be empty
occur explicitly as differences of set-valued tableaux.

\begin{lem}\label{l.inc_exc}
$\displaystyle N_S=\sum\limits_{R\subset
S}(-1)^{|\ol{R}|}-\sum\limits_{R\subset S\atop
\sstab(R)=\emptyset}(-1)^{|\ol{R}|}$.
\end{lem}
\begin{proof}
If $q_S=0$ (i.e., $\st(S)=\emptyset$), then for every set-valued
tableau $R\subset S$, $\st(R)=\emptyset$; thus the result is
trivially true. Assume $q_S\neq 0$.

For $k\geq 1$, we have that $\cN_k(S)=\cM_k(S)\setminus
\bigcup\limits_{R\subsetneqq S}\cM_k(R)$; by the
inclusion-exclusion principle,
\begin{equation*}
|\cN_k(S)|
=\sum\limits_{R\subset S}-1^{|\ol{R}|}\cdot|\cM_k(R)|
=\sum\limits_{R\subset S}(-1)^{|\ol{R}|}{|\sstab(R)| \choose k}
=\sum\limits_{R\subset S}(-1)^{|\ol{R}|}{q_R \choose k}
\end{equation*}
where we use the convention ${a\choose b}=0$ if $a<b$. Thus,
\begin{align*}
N_S&=\sum\limits_{k=1}^{q_S}(-1)^{k+1} |\cN_k(S)|\\
&=\sum\limits_{k=1}^{q_S}(-1)^{k+1}\sum\limits_{R\subset
S}(-1)^{|\ol{R}|}{q_R \choose k}\\
&=\sum\limits_{R\subset
S}(-1)^{|\ol{R}|} \sum\limits_{k=1}^{q_S}(-1)^{k+1}{q_R \choose k}.\\
\end{align*}
The result now follows from the fact that for any $R\subset S$,
$q_R\leq q_S$, and therefore
\begin{equation*}\sum\limits_{k=1}^{q_S}(-1)^{k+1}{q_R \choose
k}=
\begin{cases}
1&\text{ if }q_R\neq 0\\
0&\text{ if }q_R=0.
\end{cases}
\end{equation*}
\end{proof}

\noindent The following Lemma gives the value of the first
summation in Lemma \ref{l.inc_exc}.
\begin{lem}\label{l.inc_exc_first_term}
$\displaystyle \sum\limits_{R\subset
S}(-1)^{\ol{R}}=(-1)^{|S|+\|S\|}$.
\end{lem}
\begin{proof}
\begin{align*}
\sum\limits_{R\subset S}(-1)^{\ol{R}}
&=\sum\limits_{R_{1,1}\subset S_{1,1}\atop R_{1,1}\neq\emptyset}
\cdots \sum\limits_{R_{u,v}\subset S_{u,v}\atop
R_{u,v}\neq\emptyset}
(-1)^{|S_{1,1}\setminus R_{1,1}|\,+\,\cdots\,+\,|S_{u,v}\setminus R_{u,v}|}\\
&=\prod\limits_{i,j}\sum\limits_{R_{i,j}\subset S_{i,j}\atop
R_{i,j}\neq\emptyset}
(-1)^{|S_{i,j}\setminus R_{i,j}|}\\
&=\prod\limits_{i,j}\sum\limits_{R'_{i,j}\subsetneqq S_{i,j}}
(-1)^{|R'_{i,j}|}\\
&=\prod\limits_{i,j}\sum\limits_{k=0}^{|S_{i,j}|-1}
(-1)^{k}{|S_{i,j}|\choose k}\\
&=\prod\limits_{i,j}(-1)^{|S_{i,j}|-1}\\
&=(-1)^{|S|+\|S\|}.
\end{align*}
\end{proof}

If $S$ is semistandard, then for every $R\subset S$,
$\sstab(R)\neq\emptyset$.  Thus Lemmas \ref{l.inc_exc} and
\ref{l.inc_exc_first_term} imply Lemma \ref{l.N_S}(i).

The following Lemma gives the value of the second summation in
Lemma \ref{l.inc_exc} when $S$ is not semistandard. Lemmas
\ref{l.inc_exc}, \ref{l.inc_exc_first_term}, and the following
Lemma imply Lemma \ref{l.N_S}(ii).
\begin{lem}\label{l.inc_exc_secnd_term}
If $S$ is not semistandard, then
$\displaystyle\sum\limits_{R\subset S\atop
\sstab(R)=\emptyset}(-1)^{|\ol{R}|}=(-1)^{|S|+\|S\|}$.
\end{lem}
\begin{proof}
 Define
$\cZ(S)=\{R\subset S\mid \sstab(R)=\emptyset\}$, so that
\begin{equation}\label{e.red_1}
\sum\limits_{R\subset S\atop \sstab(R)=\emptyset}(-1)^{|\ol{R}|}=
\sum\limits_{R\in \cZ(S)}(-1)^{|\ol{R}|}.
\end{equation}
We make a series of reductions.

Let $x$ and $y$ be the row and column numbers of a box of $S$
where semistandardness is violated either on the top or left, but
not on the right or bottom. Define
\begin{align*}
\cZ'(S)&=\{R\in\cZ(S)\mid \{S_{x,y,1}\}=R_{x,y}\}\\
\cZ''(S)&=\{R\in\cZ(S)\mid \{S_{x,y,1}\}\subsetneqq R_{x,y} \}\\
\cZ'''(S)&=\{R\in\cZ(S)\mid S_{x,y,1}\not\in R_{x,y}\}.
\end{align*}
Then $\cZ(S)=\cZ'(S)\ \dot{\cup}\ \cZ''(S)\ \dot{\cup}\
\cZ'''(S)$. Consider the bijection from $\cZ''(S)$ to $\cZ'''(S)$
defined by $R\mapsto R\setminus S_{x,y,1}$.  The set-valued
tableaux paired under this bijection contribute opposite signs to
(\ref{e.red_1}). Thus
\begin{equation}\label{e.red_2}
\sum\limits_{R\in\cZ(S)}(-1)^{|\ol{R}|}= \sum\limits_{R\in
\cZ'(S)}(-1)^{|\ol{R}|}.
\end{equation}
Let $g=S_{x,y,1}$. Define
\begin{align*}
\cY'(S)&=\{R\in \cZ'(S)\mid R_{x-1,y,k}< g, R_{x,y-1,l}\leq  g,
\text{ some }k,l\}\\
\cY''(S)&=\{R\in \cZ'(S)\mid R_{x-1,y,k}\geq g, \text{ all }k\}\\
\cY'''(S)&=\{R\in \cZ'(S)\mid R_{x,y-1,k}> g, \text{ all }k\}.
\end{align*}
Then $\cZ'(S)=\cY'(S)\ \dot{\cup}\ (\cY''(S)\cup \cY'''(S))$.
Therefore
\begin{align}\label{e.red_3}
\sum\limits_{R\in \cZ'(S)}(-1)^{|\ol{R}|}&=\sum\limits_{R\in
\cY'(S)}(-1)^{|\ol{R}|}+\sum\limits_{R\in
\cY''(S)}(-1)^{|\ol{R}|}\notag\\
&\ \ \ +\sum\limits_{R\in \cY'''(S)}(-1)^{|\ol{R}|}
-\sum\limits_{R\in \cY''(S)\cap \cY'''(S)}(-1)^{|\ol{R}|}.
\end{align}
We compute the last three summations on the right hand side of
(\ref{e.red_3}).  Assume that $\cY''(S)\neq\emptyset$.  Let $S'$
be the tableaux obtained from $S$ by removing all entries other
than $g$ from $S_{x,y}$ and all entries less than $g$ from
$S_{x-1,y}$. Then $\cY''(S)=\{R\subset S'\}$. Thus
\begin{align*}
\sum\limits_{R\in \cY''(S)}(-1)^{|\ol{R}|}&=\sum\limits_{R\subset
S'}(-1)^{|S\setminus R|}\\
&=\sum\limits_{R\subset S'}(-1)^{|S\setminus S'|+|S'\setminus
R|}\\
&=(-1)^{|S\setminus S'|}\sum\limits_{R\subset
S'}(-1)^{|S'\setminus R|}\\
&=(-1)^{|S\setminus S'|}(-1)^{|S'|+\|S'\|}\\
&=(-1)^{|S|+\|S\|},
\end{align*}
where the seconds to last equality follows from Lemma
\ref{l.inc_exc_first_term}, and the last equality form the fact
that $\|S'\|=\|S\|$. In a similar manner, one shows that that
\begin{align*}
&\sum\limits_{R\in \cY'''(S)}(-1)^{|\ol{R}|}
=(-1)^{|S|+\|S\|}\ \ \hbox{ if }\cY'''(S)\neq\emptyset,\hbox{ and}\\
&\sum\limits_{R\in \cY''(S)\cap \cY'''(S)}=(-1)^{|S|+\|S\|}\
\hbox{ if }\cY''(S)\cap \cY'''(S)\neq\emptyset.
\end{align*}
Either $\cY''(S)$ or $\cY'''(S)$ must be nonempty, and if both of
them are nonempty, then so must be their intersection. It follows
that
\begin{equation}\label{e.red_5}
\sum\limits_{R\in \cZ'(S)}(-1)^{|\ol{R}|}=
(-1)^{|S|+\|S\|}+\sum\limits_{R\in \cY'(S)}(-1)^{|\ol{R}|}.
\end{equation}
Define $$\cX(S)=\{R\in \cY'(S)\mid R_{x-1,y,k}< g, R_{x,y-1,l}\leq
g, \text{ all }k,l\}.$$ Let $A=\{a_1,\ldots,a_{r}\}$ be the
entries of $S_{x-1,y}$ which are greater than or equal to $g$, let
$B=\{b_1,\ldots,b_{s}\}$ be the entries of $S_{x,y-1}$ which are
greater than $g$, and let $t=r+s$. For each $R\in\cX(S)$, define
$\cY'_R(S)$ to be all the set-valued tableaux obtained by adding
elements of $A$ to $R_{x-1,y}$ and elements of $B$ to $S_{x,y-1}$.
Then
\begin{equation*}
\cY'(S)=\dot{\bigcup\limits_{R\in\cX(S)}}\cY'_R(S),
\end{equation*}
and
\begin{align}\label{e.red_6}
\sum\limits_{R\in
\cY'(S)}(-1)^{|\ol{R}|}&=\sum\limits_{R\in\cX(S)}\sum\limits_{Q\in
\cY'_R(S)}(-1)^{|S\setminus Q|}\notag\\
&=\sum\limits_{R\in\cX(S)}(-1)^{|S\setminus
R|}\left((-1)^0{t\choose
0}+(-1)^1{t\choose 1}+\cdots+(-1)^t{t\choose t}\right)\notag\\
&=0.
\end{align}
Combining (\ref{e.red_1}), (\ref{e.red_2}), (\ref{e.red_5}), and
(\ref{e.red_6}), we obtain the result.

\end{proof}


\providecommand{\bysame}{\leavevmode\hbox
to3em{\hrulefill}\thinspace}
\providecommand{\MR}{\relax\ifhmode\unskip\space\fi MR }
\providecommand{\MRhref}[2]{%
  \href{http://www.ams.org/mathscinet-getitem?mr=#1}{#2}
} \providecommand{\href}[2]{#2}

\vspace{1em}

\noindent \textsc{Department of Mathematics, Virginia Tech,
Blacksburg, VA 24063}

\noindent \textsl{Email address}: \texttt{vkreiman@vt.edu}



\end{document}